\documentclass[11pt]{article}

\usepackage{amsmath}
\usepackage{amssymb}
\usepackage{latexsym}
\usepackage{graphicx}
\usepackage{color}
\usepackage{subfigure}

\newcommand{\eps}{\varepsilon}
\newcommand {\1}{\textbf{1}}
\newcommand{\R}{\mathcal{R}}

\title{Metastability in Stochastic Replicator Dynamics}

\author{Konstantin Avrachenkov\footnote{INRIA Sophia Antipolis,
2004 Route des Lucioles, B.P.93, 06902, Sophia Antipolis Cedex, France.
Tel: +33 4 92.38.77.51 E-mail: k.avrachenkov@inria.fr } \ and \ Vivek\ S.\ Borkar\footnote{Department
of Electrical Engineering, Indian Institute of Technology Bombay, Powai, Mumbai 400076, India. E-mail: borkar.vs@gmail.com}}

\begin{document}

\maketitle

\vspace{.5in}

\noindent \textbf{Abstract:}  We consider a novel model of stochastic replicator dynamics for potential games that converts to a Langevin equation on a sphere after a change of variables.  This is distinct from the models of stochastic replicator dynamics studied earlier. In particular, it is ill-posed due to non-uniqueness of solutions, but is amenable to a natural selection principle that picks a unique solution. The model allows us to make specific statements regarding metastable states such as small noise asymptotics for mean exit times from their domain of attraction, and quasi-stationary measures. We illustrate the general results by specializing them to replicator
dynamics on graphs and demonstrate that the numerical experiments support theoretical predictions.\\

\noindent \textbf{Key words:} stochastic replicator dynamics; Langevin equation on sphere; metastable states; intermittency; small noise asymptotics; mean exit time; quasi-stationary distributions\\

\noindent \textbf{MSC:} Primary: 91A22 Secondary: 91A25; 60H10; 34F05

\section{Introduction}

This work is motivated by the studies on intermittent dynamics
of complex large-scale systems \cite{Cetal02}, \cite{Piovani}, \cite{TY03}.
We model the dynamics of a large-scale complex system by replicator
dynamics with noise, a classical model in evolutionary biology. Replicator dynamics per se is an extensively studied dynamics \cite{Hof}, \cite{Sand}, \cite{Weibull}. Its noisy version is much less so. In their seminal work, Foster and Young \cite{Foster} considered small noise asymptotics of the stationary regime of replicator dynamics perturbed by small noise and used it to formulate the notion of \textit{stochastically stable} equilibria, i.e.,  Nash equilibria that arise as limit(s) when noise is reduced to nil in the stationary regime. This provided a selection principle for Nash equilibria (see also \cite{Young}). This is of interest because multiple Nash equilibria is a commonplace situation and one does need a selection principle of some sort to narrow down the relevant possibilities. While early work in this direction was based on static refinements of the notion of Nash equilibrium \cite{Harsanyi}, \cite{van}, later developments sought dynamic models where some Nash equilibria would naturally emerge as stable equilibria of the dynamics \cite{Fudenberg}, \cite{Sand}. This, however, only partially solves the problem, because, as mentioned above, even stable equilibria can be many. The important observation of Foster and Young was that only some of them may be stable under stochastic perturbations in the stationary regime, a physically relevant concern. (See also \cite{Manapat}, \cite{Mas} for interesting takes on behavioral aspects of noise and \cite{Borgers}, \cite{sarin}, \cite{chinese} for connections with finite population models. The general philosophy of using small noise limits for selection principles in ill-posed dynamics has been attributed to Kolmogorov in \cite{EckRu}, p.\ 626.)

The elegant machinery for small noise asymptotics of noise-perturbed dynamics due to Freidlin and Wentzell \cite{FW}  gives us a handle to tag these equilibria in a precise manner, viz., as minimizers of the so-called Freidlin-Wentzell potential associated with the process. This is precisely what Foster and Young did. But this still leaves wide open the issue of metastability, a familiar phenomenon in statistical physics \cite{Bovier}. Well known examples are supercooled water, air supersaturated with vapor, etc., which on suitable perturbation result in resp.\ freezing (ice) or condensation (rain). There have been various analytic approaches to quantitative analysis of metastability, see the opening chapters of \cite{Bovier} for a succinct overview and historical account. Our approach follows that formulated by Kramers \cite{Kramers} who viewed them as equilibria of a gradient system perturbed by noise. The later seminal work of Freidlin and Wentzell \cite{FW} generalized this viewpoint considerably, encompassing non-gradient systems and attractors more general than equilibria. In this viewpoint, metastable states  are stable equilibria for the deterministic dynamics that will not arise as a small noise limit in the stationary regime, but are still important because of the celebrated quasi-stationarity phenomenon \cite{Col}: the process may spend long time near such equilibria on time scales relevant to us before approaching true stationarity. An alternative though related viewpoint is that these are equilibria which the system gets out of when given a perturbation of a certain magnitude.  In case of gradient systems, the size of the perturbation will be related to the barrier between the current metastable state and its neighboring equilibria as captured by the optimal trajectory (i.e., one minimizing the `rate function' of the associated large deviations) needed to climb in order to get out of the current `valley' or domain of attraction for the noiseless system. Analogous ideas exist for non-gradient systems in terms of the aforementioned Freidlin-Wentzell potential. See \cite{Vares} for an extensive account of this approach and \cite{Bianchi} for some recent developments.

A related phenomenon is that of intermittency, where one sees the process spend nontrivial time in the neighborhood of a succession of stable equilibria with abrupt noise-induced episodic transitions between them (see \cite{Cetal02}, \cite{Piovani}, \cite{TY03} for examples from physics and biology). This will be the case when the noise is small but not too small so that these transitions occur on an observable time scale. This should not be confused with the phenomenon of intermittency in chaotic systems wherein a deterministic dynamics exhibits long periods of slow variation punctuated by spurts of chaotic behavior \cite{Ott}. The analysis of Fredlin-Wentzell shows that the process may be viewed as approximating a Markov chain over stable equilibria (more generally, stable attractors) with transition times and probabilities characterized in terms of a suitable `action' functional \cite{FW}. The small noise asymptotics gives a feel for the noise levels for which one expects intermittency instead of quasi-stationarity.

In view of the foregoing one important numerical quantity associated with such a metastable equilibrium is the mean exit time from its domain of attraction, or rather, the scaling properties thereof in the small noise limit. Our aim here is to present one such result in the case of replicator dynamics with potential. While bounds on this quantity are available \cite{Imhof} for the specific model considered by \cite{C00}, \cite{FudHar}, we give a finer result using the Freidlin--Wentzell theory \cite{FW}. A highlight of our approach is that it is based on an interesting connection between replicator dynamics for potential games and gradient descent on a sphere, a compact manifold without boundary which allows the application of Freidlin-Wentzell theory without facing the complications encountered by Foster and Young in doing so - see, e.g., the correction note to \cite{Foster}. It provides an alternative treatment to such issues distinct from that of \cite{Imhof} that is simpler in many ways. We also gain some important insights into quasi-stationary distributions.

In related works \cite{SandSt1}, \cite{SandSt2}, replicator dynamics is viewed as a large population `mean field' limit of finite population stochastic models (see, e.g., \cite{BorkarSund12}) and large deviations phenomena in the stationary regime are analyzed. Our emphasis, in contrast, is on quasi-stationary behavior. As noted in \cite{Vishnoi}, the quasi-stationary behavior is closely associated with the `hardness' of the underlying optimization problem of minimizing the potential. This also reflects in the slow convergence of the simulated annealing algorithm for non-convex potentials, which is a time-inhomogeneous stochastic gradient descent with slowly decreasing noise \cite{Aiba}, \cite{Catoni1}, \cite{Catoni2}. Other related works on asymptotic behavior of noisy replicator dynamics are \cite{HofImh}, \cite{Mertiko1}, \cite{Mertiko2}.  Mathematically, what sets our model apart is a diffusion matrix involving in non-diagonal positions square-roots of state variables that puts it beyond available well-posedness results for stochastic differential equations. In fact it is ill-posed due to non-uniqueness of solutions, thus we need to resort to a well motivated selection principle to pick a unique choice.

To summarize, the main contributions of the present work are: (a) we introduce a new
variant of the stochastic replicator dynamics, which in particular allows us to work around the problem of ``corners''
in \cite{Foster}, (b) we propose a transformation that converts our
model to a Langevin equation on a sphere; (c) a selection principle to overcome the ill-posedness of the original model is then proposed; (c) we characterize for our model the exit time asymptotics from  metastable states and
also characterize its quasi-stationary distributions; (d) we apply these general results to study the intermittent behaviour of large complex systems described by graphs; and finally, (e) we illustrate main theoretical findings by numerical examples.

The article is organized as follows. The next section introduces our model of stochastic replicator dynamics
and discusses in detail its relationship to other stochastic replicator dynamics in literature.
Section~3 uses a transformation to map it into the Langevin equation on a sphere and exploits the structure therein to analyze the small noise asymptotics. Section~4 specializes the results to the case of quadratic
potentials on graphs and makes some relevant observations for this special case which has important applications
such as the study of the intermittent behaviour in complex large-scale systems \cite{Cetal02}, \cite{Piovani},
\cite{TY03} and clique detection problem \cite{AKS98}, \cite{Bomze}, \cite{Bomze99}, \cite{Bomze02}.
Section~5 presents some numerical experiments to support the theoretical findings.
We conclude the paper with suggestions for future research in Section~6.\\

%
%


\newpage

\section{Stochastic replicator dynamics with potential}

Replicator dynamics, introduced by Taylor and Jonker \cite{Taylor}, is a model for behavior adjustment of population games, i.e., a game theoretic scenario in which we analyze collective behavior of a population of agents belonging to certain `types' or `species' in conflict situation, retaining only their species/type tag and erasing their individual identity. This is because their strategic behavior is dictated by the type, as also the environment they face is dictated by the collective sizes and strategic behavior of other types. The origin of this particular dynamics of population games (among many, see \cite{Sand} for an extensive treatment) lies in models of Darwinian evolution in biology \cite{Hof}. Issues of interest typically are convergence to an equilibrium behavior, persistence or non-extinction, etc. Among others, one very significant outcome of this theory is the link between the so-called evolutionarily stable states, a biologically relevant equilibrium notion introduced by John Maynard Smith, and asymptotically stable equilibrium points of replicator dynamics. See \cite{Hof}, \cite{Sand}, \cite{Weibull} for a detailed treatment. We further restrict to a special subclass of these, viz., potential games wherein the payoff function is restricted to be of a particular form, see again \cite{Sand} for an extensive account of the deterministic case.  This subclass, defined later, is much more analytically amenable and has important applications, e.g., in congestion games.

Our focus will be on \textit{stochastic} replicator dynamics with potential. We begin with a motivating example of deterministic replicator dynamics with linear payoffs and then introduce the more general paradigm.

\subsection{Deterministic replicator dynamics and potential games}

Let
$$P^n := \left\{x = [x_1, \cdots, x_n] : x_i \geq 0 \ \forall i, \ \sum_{i=1}^nx_i = 1\right\}$$
denote the $(n-1)$-dimensional probability simplex. The deterministic replicator dynamics with linear payoffs is the following differential equation in $P^n$:
\begin{equation}
\label{eq:detrepdyncont}
\frac{dx_i}{dt}(t) = x_i(t)\left([Mx(t)]_i - x^T(t)Mx(t)\right),
\end{equation}
where the $n\times n$ matrix $M$ is a payoff matrix, i.e., the $i$th element $[Mx]_i$ of $Mx$ is the payoff to species $i$ when the population profile is $x$. We assume that $M$ is symmetric.

In evolutionary biology, $x_i(t)$ is interpreted as the fraction of species $i$ in the population, $x(t)$ thus being the overall population profile.  The quantity $[Mx(t)]_i$ is the payoff for species $i$ at time $t$ when the population profile is $x(t)$. The quantity $x^T(t)Mx(t)$ is then the average payoff for the population at time $t$. Thus, the population share of species $i$ waxes or wanes at a rate proportional to its excess, resp., deficit of payoff over the population average.  The standard equilibrium notion for population game is that of Nash equilibrium defined as follows: $x^*$  is a Nash equilibrium if $x^TMx^* \leq (x^*)^TMx^*, \ \forall \ x \in P^n$.  An analogous definition is also used when $Mx$ is replaced by a nonlinear payoff function $V(x)$, as we do a bit later. A population profile $\hat{x}$ is evolutionarily stable if $\hat{x}^TMx > x^TMx$ for $x \in P^n$ sufficiently close to $\hat{x}$. Then Nash equilibria are the equilibrium points of the deterministic replicator dynamics and the evolutionarily stable states are asymptotically stable equilibria, though neither converse holds (see e.g., \cite{Hof}, \cite{Sand}, \cite{Weibull}.)  While natural and appealing, the problem with stopping with these definitions is that of multiple equilibria as already mentioned, calling for a selection principle which stochastic replicator dynamics tries to provide.\\

Equation (\ref{eq:detrepdyncont}) can be written in the following matrix-vector form
\begin{equation}
\dot{x}(t) = \mbox{diag}(x(t))(Mx(t) - x(t)^TMx(t)\1), \label{ode2}
\end{equation}
where we introduce the notation:
 \begin{itemize}
 \item diag$(z)$ for $z = [z_1, \cdots, z_n]^T \in \R^n$ is the diagonal matrix whose $i$th diagonal entry is $z_i \ \forall i;$ and,

  \item $\1$ is the constant vector of all $1$'s.
  \end{itemize}
Then, we observe that (\ref{ode2}) is a particular case of
\begin{equation}
\dot{x}(t) = \mbox{diag}(x(t))(\nabla F(x(t)) - x(t)^T\nabla F(x(t))\1), \label{ode3}
\end{equation}
with $F(x) := \frac{1}{2}x^TMx$. Given $x(t) \in P^n$, the right hand side summed over its components adds to zero, implying that $\sum_ix_i(t) \equiv$ a constant on $P^n$. Since the $i$th component of this dynamics is of the form
$$\dot{x}_i(t) = x_i(t)g_i(t)$$
for a suitably defined $g_i$, its solution is of the form $x_i(0)exp(\int_0^tg_i(s)ds) \geq 0$, which together with the preceding observation implies that  the  simplex
of probability vectors  is invariant under (\ref{ode3}). In fact, it is also clear that if a component is zero for some $t$, it is always so - this is immediate from the uniqueness of solutions.  So the faces of $P^n$ are  also invariant under (\ref{ode3}) and the trajectory in any face cannot enter a lower dimensional face in finite time.\\

For a replicator dynamics of form (\ref{ode3}), $F$ is said to be its associated potential function, assumed continuously differentiable (which it indeed is in our case). Then it is easy to check (see \cite{Weibull}, p.\ 110, or \cite{Sand}, Section~7.1, for a more general statement and an extensive account of the deterministic dynamics (\ref{ode3})) that
\begin{displaymath}
\frac{d}{dt}F(x(t)) = \sum_ix_i(t)\left(\frac{\partial F}{\partial x_i}(x(t))\right)^2 - \left( \sum_ix_i(t)\left(\frac{\partial F}{\partial x_i}(x(t))\right)\right)^2 \geq 0,
\end{displaymath}
with equality on the  set
$$
A_{eq}:= \left\{ x: \frac{\partial F}{\partial x_i} = \frac{\partial F}{\partial x_j}, \ \forall \ i \neq j,  \ i,j \in \mbox{support}(x) \right\}.
$$
This includes the critical points of $F$. Then the above inequality implies that $F$ is non-decreasing along any trajectory of (\ref{ode3}) and strictly increasing outside $A_{eq}$. This ensures convergence of all trajectories to $A_{eq}$. Note that $A_{eq}$ is specified by a set of nonlinear equations. We shall assume that it is a union of lower dimensional manifolds, which can usually be verified by using the implicit function theorem. In that case, it is clear that for any points in $A_{eq}$ other than the local maxima of $F$, there will be a direction transversal to $A_{eq}$ along which $F$ must increase (otherwise it would have been a local maximum).  (See \cite{Sand}, p.\ 273, for a related discussion.)  Thus $-F$ serves as a Liapunov function for (\ref{ode3}), implying convergence to the local maxima of $F$ for generic (i.e., belonging to an open dense set) initial conditions. Then (\ref{ode2}) is a special case with $F(x) = \frac{1}{2}x^TMx$. The entire development to follow goes through for general continuously differentiable
$F: P^n \mapsto \R$, so we shall consider this general case henceforth.

\subsection{Stochastic replicator dynamics}

The stochastic or noise-perturbed replicator dynamics is the following stochastic differential equation in $P^n$:
\begin{equation}
\label{eq:stochrepdyncont}
dx_i(t) = x_i(t)\left( \frac{\partial F}{\partial x_i}(x(t)) - x(t)^T\nabla F(x(t))\1 \right)dt + \eps x_i(t) \Gamma_i(x(t)) dW(t),
\end{equation}
where
\begin{itemize}
\item $W(t) = [W_1(t), ... , W_n(t)]^T, \ t \geq 0,$ is  the standard $n$-dimensional Wiener process;

\item  $\epsilon > 0$ is a small parameter;

\item $F(\cdot): P^n \mapsto \R$ is a continuously differentiable potential function and
$\frac{\partial F}{\partial x_i}(x)$ is a payoff to species $i$ when the population profile is $x$;

\item $\Gamma(\cdot) = [[\Gamma_{ij}(\cdot)]]_{1 \leq i, j \leq n} : P^n \mapsto \R^{n\times n}$  with $\Gamma_i(\cdot) :=$ the $i$th row thereof.
\end{itemize}
In order that
$$x(t) = [x_1(t), \cdots, x_n(t)]^T$$
 truly evolves in a probability simplex, we require that
\begin{equation}
x^T\Gamma(x) \equiv 0, \quad \forall x \in P^n. \label{boundarycondition}
\end{equation}
In \cite{Foster}, there is an additional requirement of continuity on $\Gamma(\cdot)$ which ensures well-posedness of (\ref{eq:stochrepdyncont}). It is worth noting that in the model we shall propose, we do away with this: our choice is continuous in the interior of $P^n$, but singular at its boundaries. The product $x_i\Gamma_i(x)$ is, however, continuous, but not Lipschitz at the boundary. Not surprisingly, well-posedness fails - the solutions are not unique. Note also that our state space is the simplex $P^n$ and hence the above applies only to $P^n$, the rest of $\R^n$ is irrelevant. Except for the aforementioned regularity requirement, our model does fit the Foster-Young framework \cite{Foster} in its overall structure.

In \cite{Foster}, the noise component of (\ref{eq:stochrepdyncont}) is attributed to mutation and migration. The latter is also used to justify a reflected diffusion formulation. The continuity assumption on $\Gamma$ suffices to ensure the existence of a weak solution, but since (\ref{boundarycondition}) implies in particular that the ensuing diffusion is \textit{not} uniformly nondegenerate, some additional care is needed to ensure uniqueness. For example, Lipschitz condition would suffice. That said, the variant we introduce later does away with the uniqueness requirement and opts for a selection principle instead -- more on that later.



The standard framework for analyzing diffusions described by stochastic differential equations that allows one to handle a sufficiently rich class of problems is the Stroock-Varadhan martingale formulation, which is a `weak' formulation as opposed to the original Ito formulation in terms of `strong' solutions. In this formulation, one is given a densely defined local\footnote{Extensions to non-local operators for, e.g., discontinuous processes, exist.} operator $\mathcal{A}$ on the space of bounded continuous functions on the state space, with suitable technical conditions. By a solution for a given initial condition $x$, one means a process $X(\cdot)$ defined on \textit{some} probability space so that
$$f(X(t)) - \int_0^t\mathcal{A}f(X(s))ds, \ t \geq 0,$$
is a martingale with respect to the natural filtration of $\sigma$-fields generated by $X(t), t \geq 0$, for bounded continuous $f$ on the state space in the domain of $\mathcal{A}$. The problem is well-posed if this process is unique in \textit{law}  for each $x$ and the resulting laws, tagged by $x$, satisfy the Chapman-Kolmogorov equation. In case of non-uniqueness, under mild continuity requirements one has a nonempty compact family of solutions for each $x$. A selection principle then entails picking one solution per $x$ so that the Chapman-Kolmogorov equation is satisfied. The classic work in this direction is that of \cite{Krylov}, which, however, has limitations that we shall point out later. We shall employ a different selection principle suited for our purposes, which can be developed in a more principled manner. In particular, it is in the spirit of the Kolmogorov philosophy alluded to in the introduction and thus well-motivated for game theoretic application, somewhat in the spirit of `trembling hand perfect equilibria', a refinement of Nash equilibria in static framework defined in terms of stability under random perturbations (see, e.g., \cite{FT}, \cite{Young}).

On purely mathematical grounds, an advantage of our formulation is that it is the natural extension of deterministic replicator dynamics for which the celebrated coordinate transformation of \cite{Akin}, \cite{Shah} that maps deterministic replicator dynamics for potential games to a gradient flow on sphere, extends to the stochastic case. In fact we exploit this fact in a major way in our analysis. Furthermore, it avoids artificial restrictions on the state space in the passage from deterministic to stochastic and when converted to gradient dynamics on sphere, it allows for seamless conversion from Ito to Stratonovich form and back due to absence of a correction term, allowing either formulation to be used according to convenience. All these issues are explained later in this work, as they require some preparatory material.

Note that $P^n$ is an $(n - 1)$-dimensional manifold with boundary and corners. This presented some problem
in the analysis of small noise limit for the noisy replicator dynamics in \cite{Foster},
which required the authors to use the extension of the Freidlin-Wentzell theory due to Anderson and Orey \cite{Ander} that handles such cases (see the correction note to \cite{Foster}). We propose here a simpler route by confining ourselves to an alternative  of $\Gamma$, which we motivate next.

Let $V_i(x)=\frac{\partial F}{\partial x_i}(x)$ denote the payoff accumulation rate for species $i$ when the population profile is
$x \in P^n$. Consider first the unnormalized population dynamics for the $i$th species given by
\begin{equation}
dr_i(t) = r_i(t)V_i(x(t))dt + 2\epsilon\sqrt{r_i(t)r(t)}dW_i(t), \ t \geq 0, \label{indivdynamics}
\end{equation}
where $r_i(t)$ is the population of species $i$ at time $t$ and
\begin{equation}
\label{indivdynamics1}
r(t) := \sum_{i=1}^nr_i(t), \ \ \tilde{x}_i(t) := \frac{r_i(t)}{r(t)},
\end{equation}
is resp.\ the total population at time $t$ and the population share (as a fraction of total population) of species $i$ at time $t$ ($\frac{0}{0} = 0$ by convention). We assume $V_i : P^n \mapsto \R^n$ to be Lipschitz so that well-posedness of the corresponding deterministic dynamics
$$\dot{r}(t) = r_i(t)V_i(x(t))$$
in $\{r = [r_1, \cdots, r_n]^T  : r_i \geq 0 \ \forall i; \ \sum_ir_i > 0\}$ is not an issue. We do not consider the well-posedness of (\ref{indivdynamics}) which is a degenerate diffusion with a non-Lipschitz diffusion coefficient. But we do so later in this article  for the stochastic replicator dynamics we get after normalization. Our point of departure from earlier models of stochastic replicator dynamics is the  stochastic forcing term, which is to be interpreted as follows. The noise seen by the $i$th species as a whole has incremental variance proportional to its own population and the population of each of other species (inclusive of itself, because the `infinitesimal' agents of a species also interact with other such agents in the same species), summed over the latter. In other words, it is $\propto \sum_{j=1}^nr_i(t)r_j(t)\Delta t$ where $\Delta t$ is the time increment. This is motivated by the fact that the variance of a sum of i.i.d.\ random variables scales linearly with the number thereof, therefore the net variance in the stochastic interaction of two populations should be proportional to the product of their respective sizes, assuming that each `infinitesimal' interaction between two agents picks independent noise.  This leads to the noise term $2\epsilon\sqrt{r_i(t)r(t)}dW_i(t)$ above where $\epsilon > 0$ is a scale parameter and  $W_i(\cdot)$ are independent standard Brownian motions. This is in contrast to the earlier models that add noise directly to the normalized (replicator) dynamics or to the payoffs in the unnormalized population dynamics - the noise is now attached to each (infinitesimal) interaction episode and not to any overall gross feature such as the net payoff. This makes ours a game theoretically distinct model.


A direct application of the Ito formula to (\ref{indivdynamics}) and (\ref{indivdynamics1}) leads to:\\

\noindent \textbf{Lemma 2.1} The process $\tilde{x}(\cdot) = [\tilde{x}(1), \cdots, \tilde{x}(n)]^T$ satisfies
\begin{eqnarray}
d\tilde{x}_i(t) &=& \tilde{x}_i(t)\Big(V_i(\tilde{x}(t)) - \sum_{j=1}^n\tilde{x}_j(t)V_j(\tilde{x}(t))\Big)dt \ +  \nonumber \\
&& 2\epsilon\Big(\sqrt{\tilde{x}_i(t)}dW_i(t) - \tilde{x}_i(t)\sum_{j=1}^n\sqrt{\tilde{x}_j(t)}dW_j(t)\Big) \label{normaldynamics}\\
&=& \tilde{x}_i(t)\Big[\Big(\frac{\partial F}{\partial \tilde{x}_i}(\tilde{x}(t)) - \sum_j\tilde{x}_j(t)\frac{\partial F}{\partial \tilde{x}_j}(\tilde{x}(t))\Big)dt  \nonumber \\
&& + \ 2\frac{\epsilon}{\sqrt{\tilde{x}_i(t)}}\Big(dW_i(t) - \sqrt{\tilde{x}_i(t)}\sum_j\sqrt{\tilde{x}_j(t)}dW_j(t)\Big)\Big]. \label{repdynmodified}
\end{eqnarray}

\ \\

This will be the stochastic replicator dynamics we analyze. Note that the diffusion coefficients in
$$\sqrt{\tilde{x}_i(t)}dW_i(t) - \tilde{x}_i(t)\sum_j\sqrt{\tilde{x}_j(t)}dW_j(t)$$
 are not Lipschitz near the boundary $\partial P^n$ of $P^n$, so the standard proof of well-posedness of stochastic differential equations with Lipschitz coefficients does not apply. Even the celebrated work of Watanabe and Yamada \cite{Watanabe} that  establishes well-posedness under weaker (e.g., certain H\"{o}lder continuity) conditions does not apply because of the presence of non-diagonal terms. In fact (\ref{repdynmodified}) is \textit{not} well-posed - its solution is not unique as we shall see later. We shall therefore propose a selection principle to choose a specific solution as the `natural' or `physical' solution.
For the time being, note that in the notation of (\ref{eq:stochrepdyncont}),
$$\Gamma_i(\tilde{x}(t))dW(t) = \frac{1}{\sqrt{\tilde{x}_i(t)}}\left(dW_i(t) - \sqrt{\tilde{x}_i(t)}\sum_j\sqrt{\tilde{x}_j(t)}dW_j(t)\right),$$
or,
$$\Gamma_{ij}(x) := \frac{1}{\sqrt{x_i}}(\delta_{ij} - \sqrt{x_ix_j}),$$
where $\delta_{ij}$ is the Kronecker delta. Letting $\sqrt{x} := [\sqrt{x_1}, \cdots, \sqrt{x_n}]^T$, we have $x^T\Gamma(x) = \sqrt{x}^T(I - \sqrt{x}\sqrt{x}^T) = 0$ for $x \in P^n$ as desired,
because $I - \sqrt{x}\sqrt{x}^T$ is the projection to the space $(\sqrt{x})^{\perp}$.

We take up the ill-posedness issues and subsequent analysis of this dynamics in the next section.\\

As was already mentioned, our model is a particular
instance of a general formulation in \cite{Foster}, but with different regularity conditions.
We now compare it with other models somewhat similar in spirit.
An important predecessor is the work of Fudenberg and Harris \cite{FudHar} mentioned earlier, followed by a series of works
\cite{C00,Imhof,HofImh,Khasminskii}.  In that model, they begin with a `raw' (i.e., unnormalized)  growth model for the populations with the payoffs perturbed by scaled Brownian noise,  in the spirit of the classical derivation of the deterministic replicator dynamics \cite{Weibull}. This leads to the equation (\ref{indivdynamics}) modified
by using geometric Brownian motion as follows:
$$
dr_i(t) = r_i(t) ( V_i(x(t))dt + \epsilon_i dW_i(t) ), \quad t \geq 0.
$$
When the Ito transformation is applied to obtain the s.d.e. for the population
fractions $x_i(t) := r_i(t)(\sum_jr_j(t))^{-1}$, this leads to
\begin{eqnarray*}
dx_i(t) & = & x_i(t)\Big(V_i(x(t)) - \sum_{j=1}^n x_j(t) V_j( x(t))\Big)dt \\
& & - x_i(t)\Big( \epsilon_i^2 x_i(t) - \sum_{j=1}^n \epsilon_j^2 x_j^2(t) \Big) \\
& & + x_i(t) \Big(\epsilon_i dW_i(t) - \sum_{j=1}^n x_j(t)\epsilon_j dW_j(t)\Big).
\end{eqnarray*}
Note that new terms appear in the deterministic part.

%

Recently, Mertikopoulos and Viossat \cite{Mertiko2} proposed to introduce directly the geometric-type
Brownian motion in the s.d.e.\ for the population fractions, which resulted in the following
version of the stochastic replicator dynamics:
\begin{eqnarray*}
dx_i(t) & = & x_i(t)\Big(V_i(x(t)) - \sum_{j=1}^n x_j(t) V_j( x(t))\Big)dt \\
& & + x_i(t) \Big(\sigma_i(x(t))dW_i(t) - \sum_{j=1}^n x_j\sigma_j(x(t))dW_j(t)\Big)
\end{eqnarray*}
with the diffusion coefficients $\sigma_i(x(t)), i=1,...,n$ being Lipschitz. This is a well-posed stochastic differential equation unlike ours
and while also a particular case of \cite{Foster}, it is qualitatively different in various aspects as will become apparent as we proceed.
Most  importantly, their model shares with classical models the feature that no species can become extinct (i.e., the corresponding component hit zero) in finite time, whereas ours can, and get reflected, implying fresh arrivals of the species into the system. Also, in \cite{Mertiko2} the authors address issues different from the ones we do.

\section{A transformation of the state space and analysis of metastability}

Our first objective is to overcome the difficulty in treating corners in \cite{Foster}
by using dynamics (\ref{repdynmodified}).

This is achieved through a change of variable $y(t) = \sqrt{x(t)}$. For $y \in S^n :=$ the $(n-1)$-dimensional sphere, define $\tilde{F}(y) = F([y_1^2, \cdots, y_n^2])$.
For the deterministic o.d.e.\ (\ref{ode3}),  this change of variable leads to
\begin{equation}
\dot{y}(t) = \frac{1}{4}\left(\nabla\tilde{F}(y(t)) - \langle\nabla\tilde{F}(y(t)), \ y(t)\rangle y(t)\right) = \frac{1}{4}\nabla^*\tilde{F}(y(t)), \label{ode4}
\end{equation}
where $\nabla^*$ denotes the projected gradient operator that gives gradient in the tangent space of $S^n$. The new dynamics (\ref{ode4}) is then nothing but gradient ascent on $S^n$, a compact manifold without boundary.

This transformation goes back to \cite{Akin}, \cite{Shah}, see \cite{Sand2} for a recent treatment which also discusses other related projection dynamics. The dynamics (\ref{ode4}) as defined above is well defined on the entire $S^n$, not just on its positive orthant.
Setting $x_i = y_i^2, \ \forall i,$ establishes  a one-to-one correspondence of each orthant of $S^n$ with $P^n$. Note that all local maxima, resp.\ minima of  $F$ in $P^n$ get mapped to local maxima, resp.\ minima of $\tilde{F}$ in a one to many map.  At the boundary of the orthants, the vector fields on either side will exhibit a mirror symmetry with respect to this boundary because of the nature of our transformation, in particular, the component of the driving vector field $\frac{1}{4}\nabla^*\tilde{F}(y(t))$ normal to the boundary will be symmetric across the boundary and being continuous, must vanish at the boundary. As observed earlier, the replicator dynamics in a face of $P^n$ does not reach a lower dimensional face in finite time. The same will be true for the gradient ascent above vis-a-vis boundaries of orthants, because the normal component of the above vector field vanishes at the boundary, slowing down its approach to the boundary. This, however, is not the case for the stochastic version as we shall see later.

%

Consider now the process $\tilde{y}(\cdot)$ satisfying
\begin{eqnarray}
d\tilde{y}(t) &=& \frac{1}{4}\left(\nabla\tilde{F}(\tilde{y}(t)) - \langle\nabla\tilde{F}(\tilde{y}(t)), \ \tilde{y}(t)
\rangle \tilde{y}(t)\right)dt \nonumber \\
&& + \ \epsilon (I - \tilde{y}(t)\tilde{y}(t)^T)\circ dW(t) \label{sdestrat} \\
&=& \frac{1}{4}\nabla^*\tilde{F}(\tilde{y}(t))dt + \epsilon d\widetilde{W}(t),  \label{sde2}
\end{eqnarray}
where `$\circ$' denotes the Stratonovich integral and
$$\widetilde{W}(t) := \int_0^t(I - \tilde{y}(s)\tilde{y}(s)^T)\circ dW(s), \ t \geq 0,$$
is a Brownian motion in  $S^n$. (This is a consequence of the fact that $I - yy^T$ for $y \in S^n$ is a projection to the tangent space of $S^n$ at $y$.   Note that this is an \textit{extrinsic} description of $\widetilde{W}(\cdot)$ which requires that $S^n$ be viewed as a manifold embedded in $\R^{n}$.)
The state space of $\tilde{y}(\cdot)$ is the entire $S^n$, with full support because as a diffusion in $S^n$, $\tilde{y}(\cdot)$ is non-degenerate. In particular, unlike the deterministic case, the trajectories do cross over from one orthant to another because of the Brownian component: at a point \textit{on} the boundary, the normal component of the gradient drift is zero, so the normal component of the diffusion is a pure Brownian motion, implying that it will be on either side of the boundary with probability one in an arbitrarily small time interval. \\

\noindent \textbf{Lemma 3.1} We can equivalently write (\ref{sdestrat}) as
\begin{eqnarray}
d\tilde{y}(t) &=& \frac{1}{4}\left(\nabla\tilde{F}(\tilde{y}(t)) - \langle\nabla\tilde{F}(\tilde{y}(t)), \ \tilde{y}(t)
\rangle \tilde{y}(t)\right)dt \nonumber \\
&& + \ \epsilon (I - \tilde{y}(t)\tilde{y}(t)^T)dW(t). \label{sdeito}
\end{eqnarray}
\\

\noindent \textbf{Proof:}
The Ito to Stratonovich  conversion  of a vector s.d.e.\
\begin{displaymath}
dX(t) = m(X(t))dt + \sigma(X(t))dW(t),
\end{displaymath}
written componentwise as
\begin{displaymath}
dX_i(t) = m_i(X(t))dt + \sum_j\sigma_{ij}(X(t))dW_j(t),
\end{displaymath}
goes as follows. For any twice continuously differentiable $f : \mathcal{R}^n \mapsto \mathcal{R}$, we have by Ito formula,
\begin{eqnarray*}
df(X(t)) &=& \sum_k\frac{\partial f}{\partial x_k}(X(t))\Big(m_k(X(t))dt + \sum_{\ell}\sigma_{k\ell}(X(t))dW_{\ell}(t)\Big) + \\
&& \frac{1}{2}\sum_{k,\ell,m}\sigma_{km}(X(t))\sigma_{\ell m}(X(t))\frac{\partial^2f}{\partial x_k\partial x_{\ell}}(X(t))dt
\end{eqnarray*}
and
\begin{displaymath}
d\langle f(X(\cdot)), W_{\ell}\rangle(t) \ = \ \sum_k\frac{\partial f}{\partial x_k}(X(t))\sigma_{k\ell}(X(t))dt.
\end{displaymath}
Hence
\begin{eqnarray*}
\sigma_{ij}(X(t))\circ dW_j(t) &=& \sigma_{ij}(X(t))dW_j(t) + \frac{1}{2}d\langle\sigma_{ij}(X(\cdot)), W_j\rangle(t) \\
&=& \sigma_{ij}(X(t))dW_j(t) + \frac{1}{2}\sum_k\frac{\partial \sigma_{ij}}{\partial x_k}(X(t))\sigma_{kj}(X(t))dt.
\end{eqnarray*}
Applying this to (\ref{sdestrat}), we have
\begin{equation}
(I - \tilde{y}(t)\tilde{y}(t)^T)\circ d\widetilde{W}(t) = (I - \tilde{y}(t)\tilde{y}(t)^T)d\widetilde{W}(t), \label{bm}
\end{equation}
which proves the claim. \hfill $\Box$

\ \\

\noindent \textbf{Remark:} Ito and Stratonovich\footnote{The work \cite{Khasminskii} suggests the use of Stratonovich formulation for a different reason, viz., that it correctly captures the limiting behavior of ordinary differential equations with rapidly varying random right hand side. An older work \cite{Turelli}, however, is critical of Stratonovich formulation because it misses out on some biologically relevant or realistic aspects, but as our main  motivation is not from biology, we do not address these issues.  See \cite{Hsu} for general background on diffusions on manifolds.} integrals  are certainly not always identical and there can be an $O(\epsilon^2)$ error term. The lemma states  that this term is zero for us. We emphasize this need not always be the case. The reason it works out for us is that our stochastic gradient dynamics on  sphere is driven by a pure Brownian motion  for which the diffusion coefficient is identity and therefore there is no Stratonovich correction. Our aim in stating (\ref{sdestrat}) in Stratonovich form is purely to avoid technicalities in handling the Ito version (\ref{sdeito}) in the manifold framework. The advantage of Stratonovich formulation is that it is more amenable to smooth coordinate transformations: even though we eventually work on the probability simplex, a compact manifold with boundary that admits a single local chart for the entire manifold, our passage to  it is via a stochastic differential equation on a sphere, a compact manifold \textit{without} boundary that \textit{does not} admit a single local chart.

\bigskip

We may identify a vector $[y_1, \cdots, y_n]^T$ in $S^n$ with $[|y_1|, \cdots, |y_n|]^T$ in its positive orthant. Under this equivalence relation, the quotient space is the positive orthant $S^n_+$ of $S^n$ and $\tilde{y}(\cdot)$ maps to a reflected diffusion $\check{y}(\cdot)$ in $S^n_+$. This requires that we modify (\ref{normaldynamics})-(\ref{repdynmodified}) to include a boundary term to account for the reflection, specifically a bounded variation increasing process supported on the  boundary that confines the original process to the prescribed domain. We also need to \textit{specify} a direction of reflection. This is because the Brownian and hence the diffusion trajectory is non-differentiable everywhere with probability one, therefore its direction of reflection on the boundary cannot be defined in a conventional manner. As in the case of simple reflected Brownian motion, pure reflection across a codimension one smooth surface entails normal reflection, as can be justified also by symmetry considerations: $-y_i'$ is the mirror image of $y_i'$ across the plane $y_i = 0$ and the line joining the two is normal to this plane. The other specification needed is whether the reflection is instantaneous or `sticky'. Our definition of $\check{y}(\cdot)$ implies instantaneous reflection. Both these specifications are achieved by adding to the dynamics a bounded variation process which confines the process to the desired domain while respecting the prescribed reflection direction. We describe this in detail below. This material is standard and a nice exposition can be found in \cite{Ikeda}. Thus, if $\check{y}(\cdot)$ is a reflected diffusion with normal reflection in $S^n_+$, it satisfies
\begin{equation}
d\check{y}(t) = \frac{1}{4}\nabla^*\tilde{F}(\check{y}(t))dt + \epsilon d\widetilde{W}(t) + dZ(t), \label{reflect}
\end{equation}
where $Z(\cdot) = [Z_1(\cdot), \cdots, Z_n(\cdot)]^T$ is a process (called the `local time' at the boundary) satisfying:
\begin{enumerate}
\item $Z(0) =$ the zero vector;

\item If $|Z_i|(T)$ denotes the total variation of $Z_i(t)$ on $[0, T]$, then $|Z_i|(t) < \infty \ \forall t \geq 0$ (i.e., it is a bounded variation process) and
$$|Z_i|(t) = \int_0^tI\{\check{y}_i(s) = 0\}d|Z|(s), \ t \geq 0$$
(i.e., it increases only when $\check{y}_i(t) = 0$, implying instantaneous reflection);

\item For $y \in \partial S_+^n$, let $C(y) := \{$ the inward normal $\}$ if $y$ belongs to one of the primary faces of $S^n_+$, and $:=$  the cone formed by the inward normals to the adjoining faces if $y$ belongs to the intersection of two or more primary faces of $S^n_+$. Then there exists a measurable $\gamma : [0, \infty) \mapsto \mathcal{R}^n$ such that $\gamma(t) \in C(\check{y}(t)), \ \forall \ t \geq 0,$ and
$$Z(t) = \int_0^t\gamma(\check{y}(s))d|Z|(s).$$
(This specifies the direction of reflection.)
\end{enumerate}

\noindent \textbf{Theorem 3.1} The stochastic differential equation (\ref{reflect}) has a unique strong solution.\\

\noindent \textbf{Proof:} Consider the stochastic differential equation in the convex set  $$\mathcal{C} := \{y = [y_1, \cdots, y_n] \in \R^n : \ y_i \geq 0 \ \forall i, \  \|y\|_1 \geq \eta\}$$ for some $\eta \in (0, 1)$,  given by
\begin{eqnarray}
d\breve{y}(t) &=& \frac{1}{4}\left(\nabla\tilde{F}(\breve{y}(t)) - \left\langle\nabla\tilde{F}(\breve{y}(t)), \ \frac{\breve{y}(t)}{\|\breve{y}(t)\|}
\right\rangle \frac{\breve{y}(t)}{\|\breve{y}(t)\|} \right)dt \nonumber \\
&& + \ \epsilon \left(I - \frac{\breve{y}(t)\breve{y}(t)^T}{\|\breve{y}(t)\|^2}\right)\circ dW(t) + dZ(t),  \label{reflect1.5}
\end{eqnarray}
with $Z(\cdot)$ defined analogously to the above. This is a degenerate reflected diffusion in $\mathcal{C}$ with normal reflection, with the property that it stays a.s.\ on the spherical shell $\mathcal{C} \cap\{y : \|y\| = \|\breve{y}(0)\|\}$ whenever $\|\breve{y}(0)\| > \eta$. In particular, if $\breve{y}(0) \in S^n_+$, it remains in $S^n_+$ and reduces to $\check{y}(\cdot)$. By Theorem 4.1, p.\ 175, of \cite{Tanaka}, the claim follows for (\ref{reflect1.5}) and \textit{ipso facto} for (\ref{reflect}), since $\eta > 0$ is arbitrary. \hfill $\Box$

\ \\

Recall that for an ill-posed diffusion, a selection  principle  involves picking one solution law per initial condition so that collectively they satisfy the Chapman-Kolmogorov conditions, see, e.g., the classical work of Krylov \cite{Krylov}.

\bigskip

\noindent \textbf{Theorem 3.2} Setting
\begin{equation}
\tilde{x}(t) := \check{y}(t)^2 \label{selection}
\end{equation}
gives a unique selection for a strong Markov solution to the ill-posed s.d.e.\ (\ref{normaldynamics})-(\ref{repdynmodified}).\\

\noindent \textbf{Proof:} By Ito formula,
\begin{eqnarray}
d\tilde{x}_i(t) &=& \tilde{x}_i(t)\Big(V_i(\tilde{x}(t)) - \sum_{j=1}^n\tilde{x}_j(t)V_j(\tilde{x}(t))\Big)dt   \nonumber \\
&& + \ 2\epsilon\Big(\sqrt{\tilde{x}_i(t)}dW_i(t) - \tilde{x}_i(t)\sum_{j=1}^n\sqrt{\tilde{x}_j(t)}\Big)dW_j(t) \nonumber \\
&& + \ 2\sqrt{\tilde{x}_i(t)}dZ_i(t),  \nonumber \\
&=& \tilde{x}_i(t)\Big(V_i(\tilde{x}(t)) - \sum_{j=1}^n\tilde{x}_j(t)V_j(\tilde{x}(t))\Big)dt   \nonumber \\
&& + \ 2\epsilon\Big(\sqrt{\tilde{x}_i(t)}dW_i(t) - \tilde{x}_i(t)\sum_{j=1}^n\sqrt{\tilde{x}_j(t)}\Big)dW_j(t)   \nonumber  \\
&=& \tilde{x}_i(t)\Big[\Big(\frac{\partial F}{\partial \tilde{x}_i}(\tilde{x}(t)) - \sum_j\tilde{x}_j(t)\frac{\partial F}{\partial \tilde{x}_j}(\tilde{x}(t))\Big)dt  \nonumber \\
&& + \ 2\frac{\epsilon}{\sqrt{\tilde{x}_i(t)}}\Big(dW_i(t) - \sqrt{\tilde{x}_i(t)}\sum_j\sqrt{\tilde{x}_j(t)}dW_j(t)\Big)\Big],   \nonumber
\end{eqnarray}
for $t \geq 0$, where we have used the fact that $2\sqrt{\tilde{x}_i(t)}dZ_i(t) \equiv 0$ because $Z_i(t)$ increases only on the set $\{t : \tilde{x}_i(t) = 0\}$. Thus $\tilde{x}(\cdot)$ defined by (\ref{selection})  satisfies (\ref{normaldynamics}) and (\ref{repdynmodified}). Furthermore, by Theorem 3.1, (\ref{selection}) uniquely specifies $\tilde{x}(\cdot)$. On the other hand, $\tilde{y}(\cdot)$ hits the boundary of $S^+_n$ and therefore $\tilde{x}(\cdot)$ must hit the boundary of $P^n$. Then any solution that remains confined to a face of $P^n$ after hitting the same with some prescribed probability $> 0$ is also a solution. This shows that the above $\tilde{x}(\cdot)$ is not a unique solution. In other words, (\ref{normaldynamics}) and (\ref{repdynmodified}) are ill-posed. Hence (\ref{selection}) implies a genuine selection process.  The  strong Markov property follows from that of the well-posed reflected diffusion $\check{y}(\cdot)$ because the two are interconvertible by a continuous invertible transformation $S^n_+ \leftrightarrow P_n$. \hfill $\Box$

 \ \\

\noindent \textbf{Remarks} 1.  This interconvertibility also implies that $\tilde{x}(\cdot)$ is adapted to the increasing filtration of  $\sigma$-fields generated by the driving Brownian motion (because $\tilde{y}(\cdot)$ is) and since it has continuous paths, is also predictable with respect to it. Thus it remains amenable to the use of standard tools of Ito calculus such as the Ito formula, a fact we shall use implicitly throughout what follows.

2. Observe that by
 the well-posedness of (\ref{reflect}) in the strong and therefore weak sense, $\check{y}(\cdot)$ is also the unique limit in law as  $0 < \eta \downarrow 0$ of the well-posed non-degenerate reflected diffusion $\check{y}^{\eta}(\cdot)$ in $(\R^n)^+$ given by
\begin{displaymath}
d\check{y}^{\eta}(t) = \frac{1}{4}\nabla^*\tilde{F}(\check{y}^{\eta}(t))dt + \epsilon d\widetilde{W}(t) + \eta d\widehat{W}(t) + dZ(t),
\end{displaymath}
where $\widehat{W}(\cdot)$ is an independent Brownian motion in $\R^n$ and $\eta > 0$. Setting $\tilde{x}_i^{\eta}(t) = \check{y}_i^{\eta}(t)^2 \ \forall i$, $\tilde{x}^{\eta}(\cdot)$ satisfies the stochastic differential equation in $(\R^n)^+$ given by
\begin{eqnarray*}
d\tilde{x}^{\eta}_i(t)
&=& \tilde{x}^{\eta}_i(t)\Big(V_i(\tilde{x}^{\eta}(t)) - \sum_{j=1}^n\tilde{x}^{\eta}_j(t)V_j(\tilde{x}^{\eta}(t))\Big)dt + \ 2\Big(\epsilon\sqrt{\tilde{x}^{\eta}_i(t)}dW_i(t) + \\
&&\eta\sqrt{\tilde{x}^{\eta}_i(t)}d\widehat{W}_i(t) - \epsilon\tilde{x}^{\eta}_i(t)\sum_{j=1}^n\sqrt{\tilde{x}^{\eta}_j(t)}dW_j(t)\Big) + \eta^2dt.
\end{eqnarray*}
This is indeed well posed as an $(\R^n)^+$-valued diffusion in the sense of Stroock-Varadhan martingale problem, in particular it has a unique weak solution \cite{BassPerkins}. Furthermore, we recover $\tilde{x}(\cdot)$ as the unique limit in law thereof as $\eta\downarrow 0$ (because the corresponding limit in law of $\check{y}^{\eta}(\cdot)$ is unique).  This selection is along the lines of the selection principle enunciated in \cite{selection}, but cannot be directly deduced from the results of \cite{selection} because (\ref{normaldynamics}) and (\ref{repdynmodified}) do not satisfy the conditions stipulated therein. Note also that the selection procedure employed above is in the spirit of the Kolmogorov philosophy of `selection through small noise limit' alluded to in the introduction and leads to a unique choice, unlike the procedure of \cite{Krylov} where the result depends on the choice  and ordering of the `test functionals' employed. It is also worth noting that a key condition, condition (1.4) of \cite{BassPerkins}, fails here for $\eta = 0$ and sure enough, there is no uniqueness as already observed.

\bigskip

We also lose the phenomenon
of extinction in this framework. In fact, our model describes a community of fixed size where different sub-communities wax or wane in their relative share depending on the selection mechanism, but do not become extinct - they can bounce back from zero population.


\bigskip

We next consider the asymptotic behavior of the above mentioned processes.

\bigskip

\noindent \textbf{Theorem 3.3} The process $\{\tilde{y}(t), t \in [0,\infty) \}$ is ergodic with stationary distribution
\begin{equation}
\xi^\epsilon(dy) = \varphi^{\epsilon}(y)\kappa(dy) := Z^{-1}e^{\frac{\tilde{F}(y)}{8\epsilon^2}}\kappa(dy) \label{statdistr}
\end{equation}
for $\kappa :=$ the uniform surface measure on $S^n$, where $Z$ is the normalizing factor
$$
Z := \int e^{\frac{\tilde{F}(y')}{8\epsilon^2}}\kappa(dy').
$$
\noindent \textbf{Proof:} Note that $\xi^\epsilon(dy)$ is the Gibbs distribution with potential $-\tilde{F}$. Its invariance is easily established by verifying that the stationary Fokker-Planck equation $\mathcal{L}^*\varphi \equiv 0$, where
$$
\mathcal{L} := \frac{1}{4}\nabla\widetilde{F}\cdot\nabla + \frac{\epsilon^2}{2}\Delta_S
$$
is the generator of the diffusion $\tilde{y}(\cdot)$ and $\mathcal{L}^*$ is its formal adjoint, given by
$$
\mathcal{L}^*(f) := -\frac{1}{4}\nabla(f\nabla \widetilde{F}\cdot \textbf{1}) + \frac{\epsilon^2}{2}\Delta_Sf
$$
for $f \in C_{\infty}(S^n)$, $\textbf{1} := [1,1,\cdots,1]^T$. ($\Delta_S$ denotes the Laplace-Beltrami operator on $S^n$.) Ergodicity follows by classical arguments: As a non-degenerate diffusion on $S^n$, $\tilde{y}(\cdot)$  has for each $t > 0$ a transition probability that is absolutely continuous with respect to $\kappa(dy)$ with strictly positive density (say) $p(y|x,t)$. Then every invariant distribution $\eta(dy)$ must satisfy
$$\int p(y|z, t)\eta(dz)\kappa(dy) = \eta(dy)$$
by invariance. It follows that any two candidate $\eta(\cdot)$ are mutually absolutely continuous w.r.t.\ $\kappa$ and therefore w.r.t.\  each other. From the Doeblin decomposition from ergodic theory of Markov processes, we know that extremal invariant measures are singular with respect to each other. The foregoing then implies that there is only one invariant measure. \hfill $\Box$

\ \\

For an extensive treatment of invariant measures of reflected diffusions, see \cite{Kang}. In particular, we have\\

\noindent \textbf{Corollary 3.1} The process $\tilde{x}(\cdot)$  has a stationary distribution supported on the interior of $P^n$ given by
\begin{equation}
\hat{Z}^{-1}\left(\prod_{m=1}^n\frac{1}{\sqrt{x_m}}\right)e^{\frac{F(x)}{8\epsilon^2}}\zeta(dx), \label{statdistr1}
\end{equation}
with $\zeta :=$ the normalized Lebesgue measure on $P^n$, $\hat{Z}$ being the normalizing factor. In particular, as $\epsilon\downarrow 0$,  the stationary distribution concentrates on the global maxima of $F$, rendering them stochastically stable.\\

\noindent \textbf{Proof:} The first claim follows by an application of the change of variables formula  to (\ref{statdistr}). The second claim then follows by standard arguments. \hfill $\Box$

\ \\

Note that we proved ergodicity for the process $\tilde{y}(\cdot)$. That of $\tilde{x}(\cdot)$ follows. The latter, however, is a statement regarding the specific selection. Other solutions remain, e.g., those restricted to the boundary, but we do not consider them. An important advantage of the explicit Gibbs distribution above is that it lays bare the metastable states as local maxima of  $F$ and flags the global maxima as stochastically stable states.

We analyze next the important phenomenon of metastability. The advantage of taking the present route is twofold. First, working with a compact manifold without boundary, we are spared the technicalities due to a non-smooth boundary which occurred in \cite{Foster}. Secondly, as in \cite{Foster},
we shall be seeking small noise asymptotics for
 concentration on the minima of the associated Freidlin-Wentzell potential (\cite{FW}, Chapter 6), which for the above gradient ascent, is simply proportional to $-F$. Hence for this special case, one can get a good handle on exit times from the domain of attraction of a stable equilibrium $x^*$ of (\ref{ode3}):

\bigskip

\noindent \textbf{Lemma 3.2} Let $\tilde{G}$ be the domain of attraction of a stable equilibrium $y^* \in S^n$ of (\ref{ode4}) with boundary $\partial \tilde{G}$ and let
$$
\tilde{\tau}_{\epsilon} := \inf\{t \geq 0 : \tilde{y}(t) \notin \tilde{G}\},
$$
where $\tilde{y}(t)$ is given by o.d.e. (\ref{sdeito}).
Then for $y \in \tilde{G}$,
\begin{equation}
\lim_{\epsilon\downarrow 0}\epsilon^2\log E_y\left[\tilde{\tau}_{\epsilon}\right]
= \frac{1}{2}\min_{z \in \partial \tilde{G}}(\tilde{F}(y^*) - \tilde{F}(z)).
\label{exitrate}
\end{equation}

\noindent \textbf{Proof:} Let $\tilde{G}_0 := \{y \in \tilde{G} : \tilde{F}(y) \geq \tilde{F}_0 := \max_{y' \in \partial \tilde{G}}\tilde{F}(y')\}$. Then it is easy to see that $\sup_{\tilde{G}}\tilde{F} = \sup_{\tilde{G}_0}\tilde{F}$. By combining Theorem 3.1, p.\ 100, and Theorem 4.1, p.\ 106,
of \cite{FW}, we then have
$$
\lim_{\epsilon\downarrow 0}\epsilon^2\log E_x\left[\tilde{\tau}_{\epsilon}\right] = \frac{1}{2}\min_{y \in \partial \tilde{G}}(\tilde{F}(y^*) - \tilde{F}(y)).
$$
The claim follows.  \hfill $\Box$

\ \\

\noindent \textbf{Theorem 3.4} If $G$ is the domain of attraction of a stable equilibrium  $x^*$ (possibly on the boundary) of (\ref{ode3}) and let
$$
\tau_{\epsilon} := \inf\{t \geq 0 : \tilde{x}(t) \notin G\}.
$$
where $\tilde{x}(t)$ is given by (\ref{selection}).
Then for $x \in G$,
\begin{equation}
\lim_{\epsilon\downarrow 0}\epsilon^2\log E_x\left[\tau_{\epsilon}\right]
= \frac{1}{2}\min_{z \in \partial G}(F(x^*) - F(z)).
\label{exitrate}
\end{equation}

\noindent \textbf{Proof:} In Lemma 3.2, let $y^* = \sqrt{x^*}$.
The claim then follows from the relationship between $F$ and $\tilde{F}$.  \hfill $\Box$

\bigskip

It should be emphasized that this captures the dominant part of the small noise asymptotics for mean exit times. A finer analysis is possible by combining Theorem 11.2, p.\ 267; Proposition 11.7, p.\ 275, and Theorem 11.12, p.\ 280, of \cite{Bovier}.

Note also that we are able to make a statement about boundary equilibria for which a general theory appears unavailable. We are able to do so because we first derive our results for a process on a manifold without boundary and then deduce the result for the desired process by a change of variables.

Following \cite{FW}, Chapter 6, we can infer considerable additional information regarding the behavior of the process for small $\epsilon$. For example, we know that with probability approaching $1$ as $\epsilon\downarrow 0$, the exit from the domain of attraction of $x^*$ occurs near a minimizer of the r.h.s.\ in (\ref{exitrate}). This allows us to consider possible paths from $x^*$ to a `ground state' or global minimum of the potential via a succession of stable equilibria, subject to the transition between two consecutive ones being compatible with the foregoing. By Theorem 4.2, Chapter 4 of \cite{FW}, we also know that the exit time distribution for each successive jump is approximately exponential. Furthermore, the process behaves like a continuous time Markov chain on the set of metastable states, i.e., the stable equilibria. Then each path is traversed with time approximately distributed as sums of independent (by the strong Markov property) approximately exponential random variables, of which the smallest exponential rate will dominate in the $\epsilon\downarrow 0$ limit. Thus the transit time from $x^*$ to a ground state will be approximately exponential with exponential rate equal to the largest, over the aforementioned paths, of the smallest exponential rate along the path. This is akin to the considerations that characterize the optimal constant for logarithmic cooling in simulated annealing \cite{Hajek}.

We conclude this section with a remark on quasi-stationary distributions. In the small noise regime, the process may spend a long time in the domain of attraction $G$ of a stable equilibrium, and thereby be well approximated by another process confined to $\bar{G}$. The stationary distribution of the latter then well approximates the behavior of the original process in the time scales of interest. In discrete Markov chain case, several notions of quasi-stationary distributions exist, but they agree in the small noise limit \cite{ABN}. In case of diffusions, one of these notions is the most convenient to work with: the limit in law as $T\uparrow\infty$ of the conditional law of the process conditioned on never exiting the region of interest till $T$. This is the so called `Q-process' of \cite{Col}, Chapter 7. As shown in \cite{Pinsky}, this turns out to be the stationary distribution of another diffusion which lives in $\bar{G}$ and whose extended generator can be exactly characterized, viz., it is  given as the unique minimizer of the (convex, lower semi-continuous) Donsker-Varadhan rate function for large deviations of empirical measures of the original diffusion from \textit{its} stationary distribution, when the minimization is performed over probability measures supported in $\bar{G}$. The results of \cite{Pinsky} are for diffusions in $\R^d$, but similar results will hold for diffusions on a sphere.
A  recent work \cite{Champ} provides an alternative route to the Q-process for diffusions which is more explicit. Translated into our framework, their Theorems 3.1-3.2 state that the extended generator of the Q-process is given by
$$\breve{\mathcal{L}}f =  \frac{\mathcal{L}(\phi f)}{\phi} + \lambda_0f,$$
where $\mathcal{L}$ is the extended generator of the original diffusion (in our case, (\ref{sde2})) and $\phi, \lambda_0$ are respectively the principal eigenfunction and eigenvalue associated with $-\mathcal{L}$ with Dirichlet boundary condition on $\partial G$. In particular, a simple application of Girsanov transformation shows that this amounts to changing the drift of (\ref{sde2}) by an additional term
$$\epsilon^2\nabla^* \log\phi.$$
We summarize this below as a theorem.\\

\noindent \textbf{Theorem 3.5} The law of the Q-process corresponds to the diffusion  in $G$ given by
\begin{displaymath}
d\tilde{y}(t)
= \frac{1}{4}\nabla^*\tilde{F}(\tilde{y}(t))dt + \epsilon^2\nabla^* \log\phi(\tilde{y}(t)) + \epsilon d\widetilde{W}(t),
\end{displaymath}
where $\phi$ is the principal eigenfunction of the operator $-\mathcal{L}$ above with Dirichlet boundary condition on $\partial G$.

\section{Intermittent dynamics on large graphs}

We next describe our general results for the case of quadratic potential $F(x) := \frac{1}{2}x^TMx$
with $x$ belonging to a high dimensional space and $M$ representing a network structure.
We know at least two important applications motivating such a setting. The first motivation comes
from the Tangled Nature model \cite{Cetal02} or its high level abstraction based on replicator
dynamics \cite{Piovani,TY03}. These models have been proposed to explain the intermittent behaviour
of complex ecological systems as well as sudden, fundamental changes in economic systems.
The authors of \cite{Piovani,TY03} considered a discretized version of the standard replicator
dynamics with matrix $M$ representing a large weighted random graph. They add to the standard
replicator dynamics the processes of extinction and mutation. The extinction happens when the
population density of a species drops below a certain threshold. Also, at each time step,
a fraction of a species can mutate to some other species with a certain probability distribution,
which can be tuned to encourage or discourage mutation to similar species. We feel that our
continuous time model share important features with the models of \cite{Cetal02,Piovani,TY03}.
As in \cite{Cetal02,Piovani,TY03}, we too observe in our model an intermittent behaviour
when a state stable for relatively long time changes suddenly to a new state, which can be
quite different in terms of its set of strategies that have a significant population share relative to the noise level. As in
\cite{Cetal02,Piovani,TY03}, our model is not confined to the interior of the simplex.
Also, as in the above cited models, the metastable states in our model have only a small
number of strategies with significant population share. The reader will be able to see these features in the
numerical examples described in the ensuing section.

Our second motivating application comes from the series of works \cite{Bomze,Bomze99,Bomze02}
on the application of the replicator dynamics to Maximum Clique Detection Problem (MCP). The goal
of MCP is to find in a graph a clique, i.e., a complete subgraph, with maximum size.
Let $A$ be the adjacency matrix of a graph $G=(V,E)$, with $V$ as the vertex set, $|V|=n$, and $E$ as the edge set, $|E|=m$. In other words, $A$ is an $n\times n$ symmetric matrix whose $(i,j)$th element $= 1$
if $(i,j) \in E$ and $= 0$, otherwise. We assume that there are no self-loops,
i.e., $A_{ii} = 0, \ \forall i$.  If we take $M=A+\frac{1}{2}I$,
the strict local maxima of $F(x)$ over the simplex and therefore Evolutionarily Stable Strategies (ESS),
will correspond to maximal cliques, see e.g. \cite{Bomze,Bomze99,Bomze02}.
Specifically, if $C$ is a maximal clique (i.e., not a proper subset of another clique) and
$x_C$ is an associated characteristic vector with mass $1/|C|$ at each vertex
belonging to $C$ and with zero mass at all other vertices, then $F$ defined over the simplex
achieves a local maximum at $x_C$. Clearly, $F$ achieves the global maximum over the simplex
at the characteristic vector corresponding to the maximum clique (the largest clique in the graph.)\\

Note that maximal cliques will typically correspond to boundary equilibria, for which our model is particularly suited, because they become  interior equilibria when the dynamics is lifted to the sphere. In general we expect the
analysis of such equilibria to be much harder in models where the boundary is not reached in finite time.\\

Now let us provide an upper bound on the expected intermittence times by specializing further
the results of Theorem 3.4. Let $x^* = x_C$ in Theorem 3.4 with $G$ its domain of attraction and $\tau_{\epsilon}$ defined correspondingly. \\

\noindent \textbf{Theorem 4.1} For any graph we have the following bound on the asymptotics of the
expected exit times from metastable states
$$
\lim_{\epsilon\downarrow 0}\epsilon^2\log E_x\left[\tau_{\epsilon}\right] \le \frac{1}{2}\left[\left(1-\frac{1}{2|C|}\right)-\frac{1}{4n}\right].
$$

\ \\

\noindent \textbf{Proof:}  Recall from \cite{Bomze08} the following
lower bound for the quadratic function $F(x)$ over the simplex:
$$
F(x) \ge \frac{1}{2}\left(\sum_{i=1}^n m_{ii}^{-1}\right)^{-1}.
$$
Since in our case $m_{ii}=1/2$, we obtain
$$
F(x) \ge \frac{1}{4n}.
$$
Using the above bound, we can write
\begin{equation}
\label{Texitbound}
\lim_{\epsilon\downarrow 0}\epsilon^2\log E_x\left[\tau_{\epsilon}\right]
=\frac{1}{2}\min_{y \in \partial G}(F(x^*) - F(y))
\le \frac{1}{2}\left[\left(1-\frac{1}{2|C|}\right)-\frac{1}{4n}\right],
\end{equation}
where $C$ is the clique associated with a local maximum.
The claim follows. \hfill $\Box$

\ \\

This suggests that the maximum mean exit time from the domain of attraction of the  local maxima
grows very slowly with $n$, implying that while
the number of local maxima increases, the sizes of the barrier heights for their domains of attraction
do not increase drastically.


We can make the result even more precise if we consider a particular model of large graphs,
e.g., the Erd\H{o}s-R\'{e}nyi $G(n,p)$ random graph model \cite{Bollobas}. Recall that in $G(n,p)$ random graph
of size $n$, the edges between any two vertices are drawn at random with probability $p$.\\

\noindent \textbf{Theorem 4.2} For sufficiently large $n$, with high probability we have
$$
\lim_{\epsilon\downarrow 0}\epsilon^2\log E_x\left[\tau_{\epsilon}\right]
=\frac{1}{2}\min_{y \in \partial G}(F(x^*) - F(y))
\le \frac{1}{2}\left[\left(1-\frac{1}{2\lceil2\log_{1/p}(n)\rceil}\right)-\frac{1}{4n}\right].
$$

\ \\

\noindent \textbf{Proof:} It is known \cite[Chapter~11]{Bollobas} that in $G(n,p)$ the maximum clique
size is either $\lfloor2\log_{1/p}(n)\rfloor$ or $\lceil2\log_{1/p}(n)\rceil$, with high probability.
By high probability we mean that when $n$ goes to infinity,
the probability of an event goes to one. The claim then is a simple rephrasing of (\ref{Texitbound})
in view of preceding remarks.
\hfill $\Box$

\section{Numerical examples}

Before proceeding to the description of two particular numerical examples, let us describe the overall discretization scheme that is
used for numerical experiments:
\begin{eqnarray*}
\widehat{Y}((n+1)\Delta) &=& Y(n\Delta) \\
&& - \ \frac{\Delta}{4} \Big\{\nabla \tilde{F}(Y(n\Delta)) - \langle Y(n\Delta), \nabla \tilde{F}(Y(n\Delta))\rangle Y(n\Delta)\Big\} \\
&&  + \  \epsilon\sqrt{\Delta}(I - Y(n\Delta)Y(n\Delta)^T)(W((n+1)\Delta) - W(n\Delta)), \\
Y_i((n+1)\Delta) &=& \frac{\widehat{Y}_i((n+1)\Delta)}{\|\widehat{Y}((n+1)\Delta)\|}, \\
X_i((n+1)\Delta) &=& Y_i^2((n+1)\Delta).
\end{eqnarray*}
Here the components of $W((n+1)\Delta) - W(n\Delta)$ are i.i.d.\ $N(0,1)$ for all $n \ge 0$.
The second step above is a renormalization step to ensure that the iterates remain on the unit sphere despite
numerical approximations. This is an example of a \textit{retraction}, a common operation for, e.g., optimization algorithms on an embedded manifold (see, e.g., \cite{Absil}, section 4.1).

\begin{figure}[h]
\begin{center}
\subfigure[The fraction of the first subpopulation, $X_1(t)$.]{\includegraphics[width=\linewidth,height=0.2\textheight]{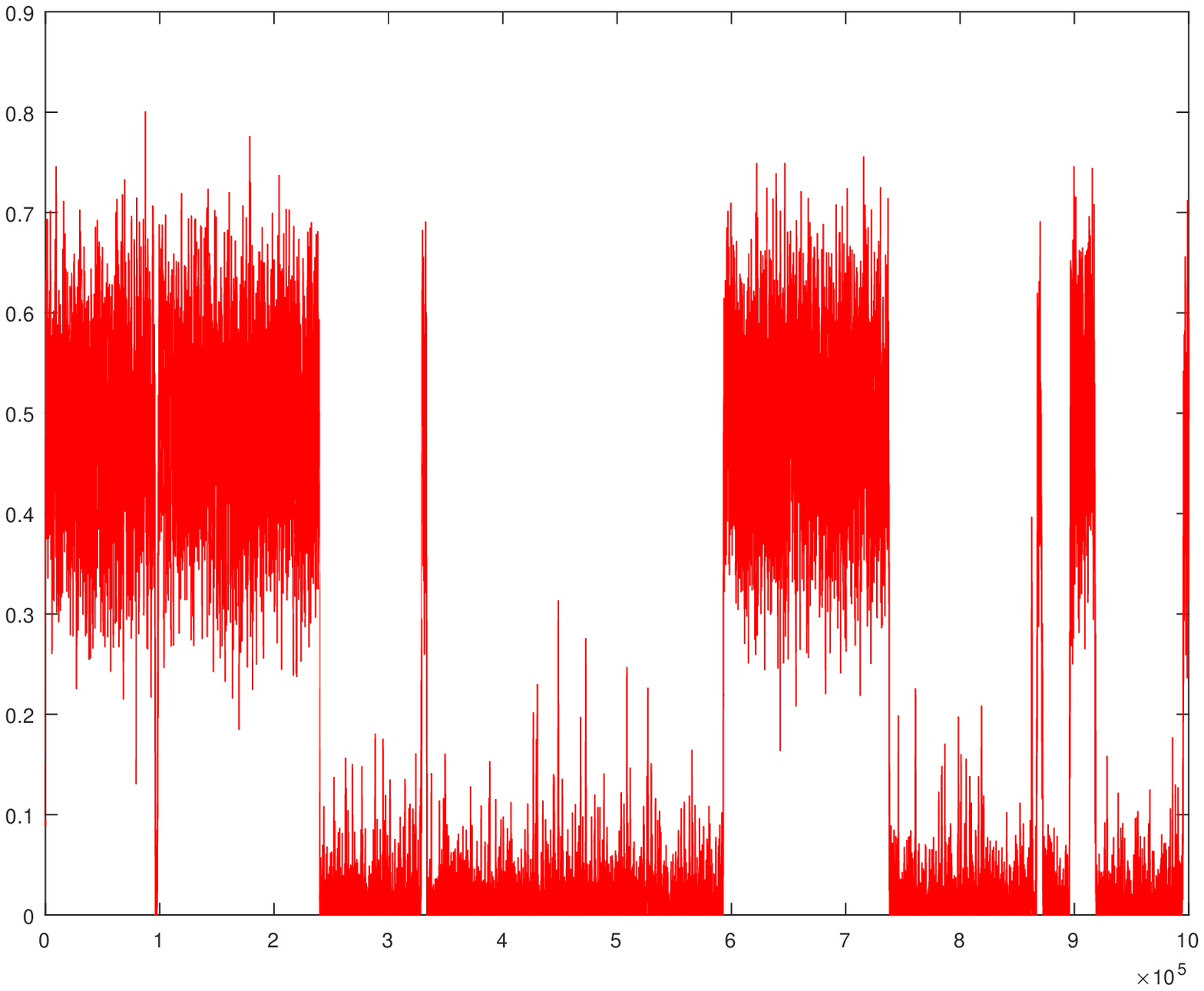}}\\
\subfigure[The fraction of the second subpopulation, $X_2(t)$.]{\includegraphics[width=\linewidth,height=0.2\textheight]{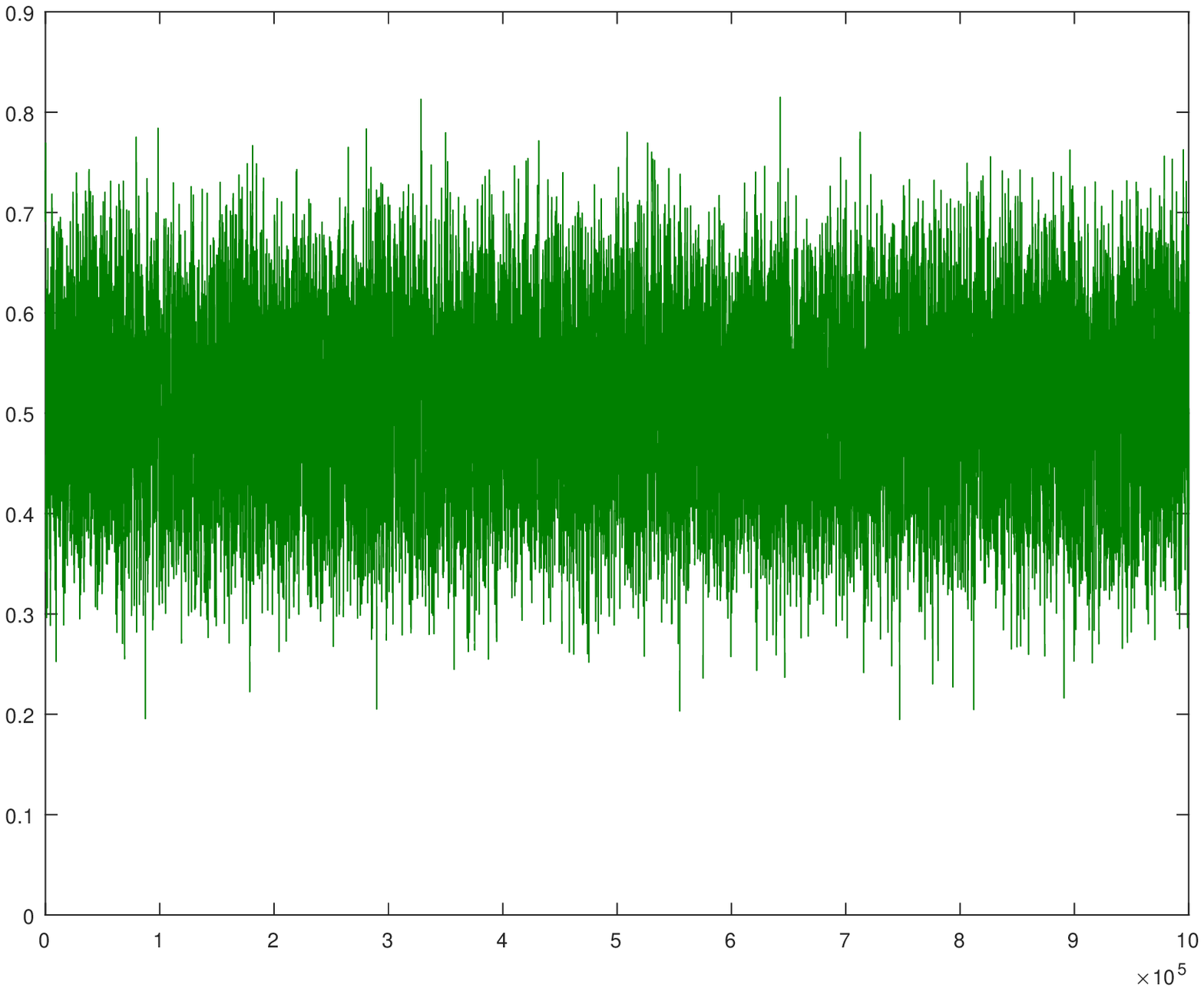}}\\
\subfigure[The fraction of the third subpopulation, $X_3(t)$.]{\includegraphics[width=\linewidth,height=0.2\textheight]{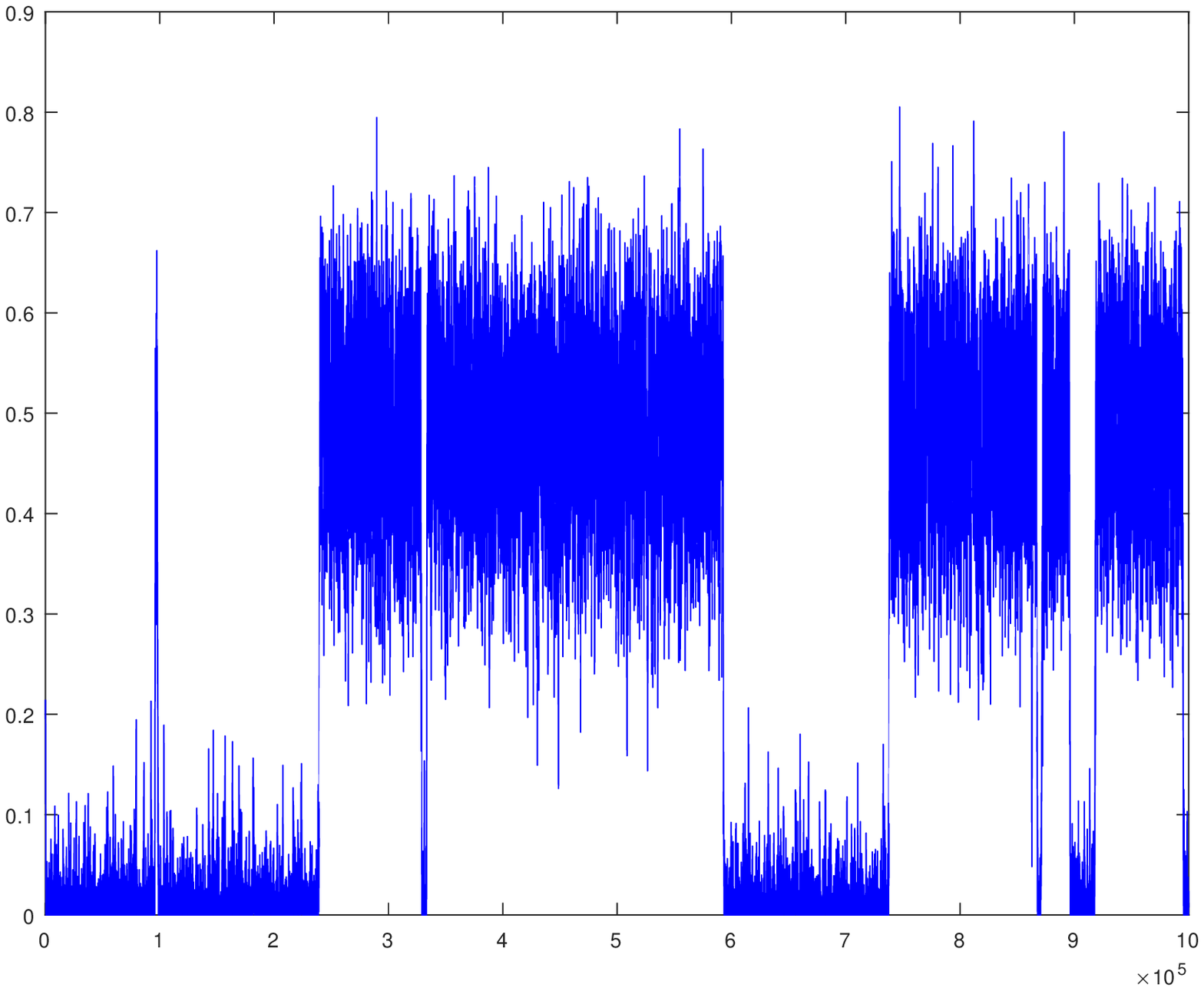}}\\
\caption{Sample of trajectory in two edge example. The x-axis
corresponds to the time and the y-axis to the subpopulation fraction.}
\label{fig:A3e005Traj}
\end{center}
\end{figure}

\begin{figure}[h]
\begin{center}
\includegraphics[scale=0.3]{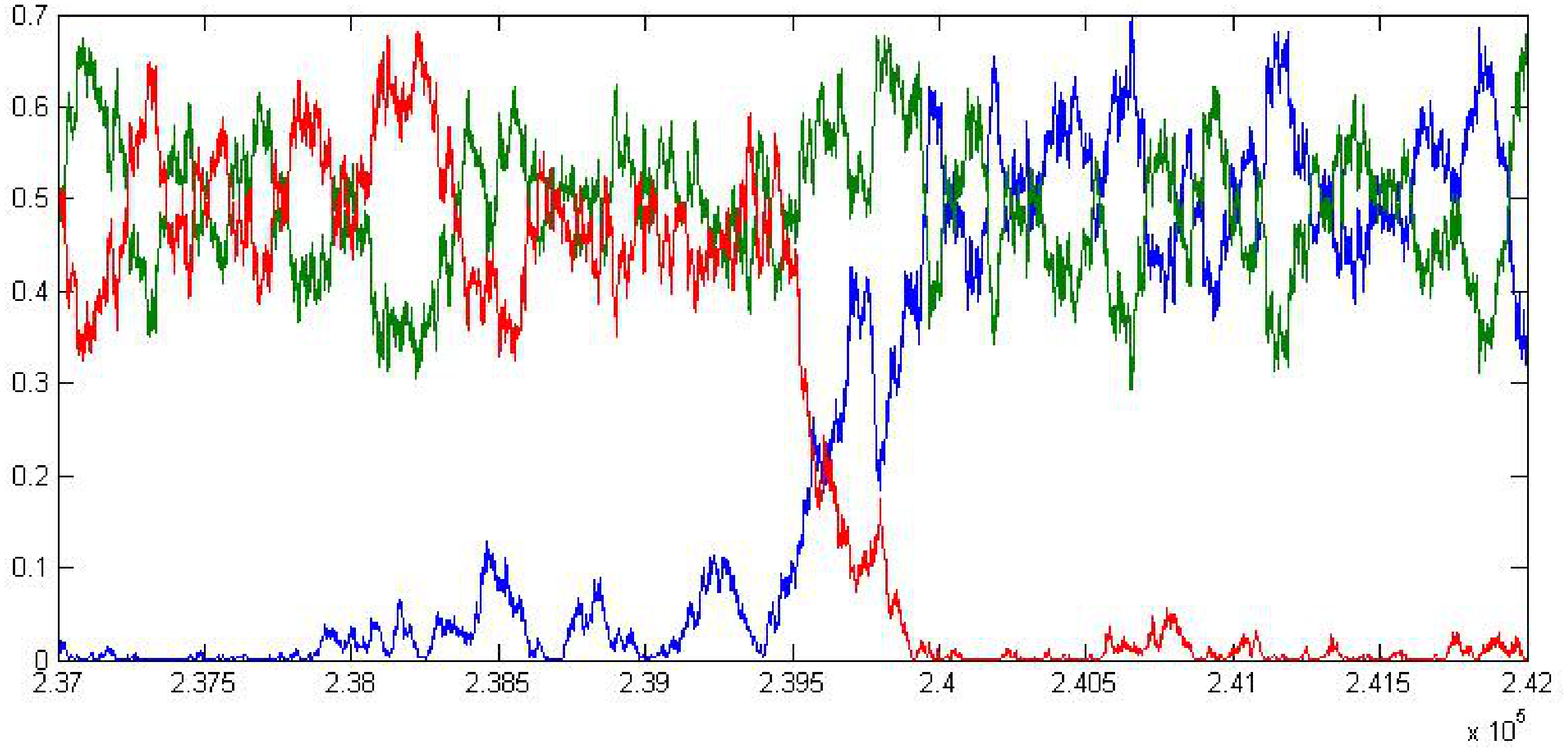}
\caption{Sample of trajectory in two edge example (zoom). The x-axis
corresponds to the time and the y-axis to the subpopulation fractions.}
\label{fig:A3e005TrajZoom}
\end{center}
\end{figure}

Let us first consider the simplest non-trivial example with two edges connected in one node.
Equivalently, it is described by the following interaction matrix:
$$
M = \left[\begin{array}{ccc}
1/2 & 1 & 0\\
1 & 1/2 & 1\\
0 & 1 & 1/2
\end{array}\right].
$$
An example of system trajectory is presented in Figure~\ref{fig:A3e005Traj} with a zoom of one
of transitions between the metastable states shown in Figure~\ref{fig:A3e005TrajZoom}.


The associated potential function $F(x)=\frac{1}{2} x^T M x$ has two global maxima
at $[ 1/2 \ 1/2 \ 0 ]^T$ and $[ 0 \ 1/2 \ 1/2 ]^T$, corresponding to the two maximum cliques.
In game theoretic terminology, these are two evolutionarily stable strategies.
In this simple example, two maximal, and also maximum, cliques are just the two edges.

\begin{figure}[h]
\begin{center}
\includegraphics[scale=0.35]{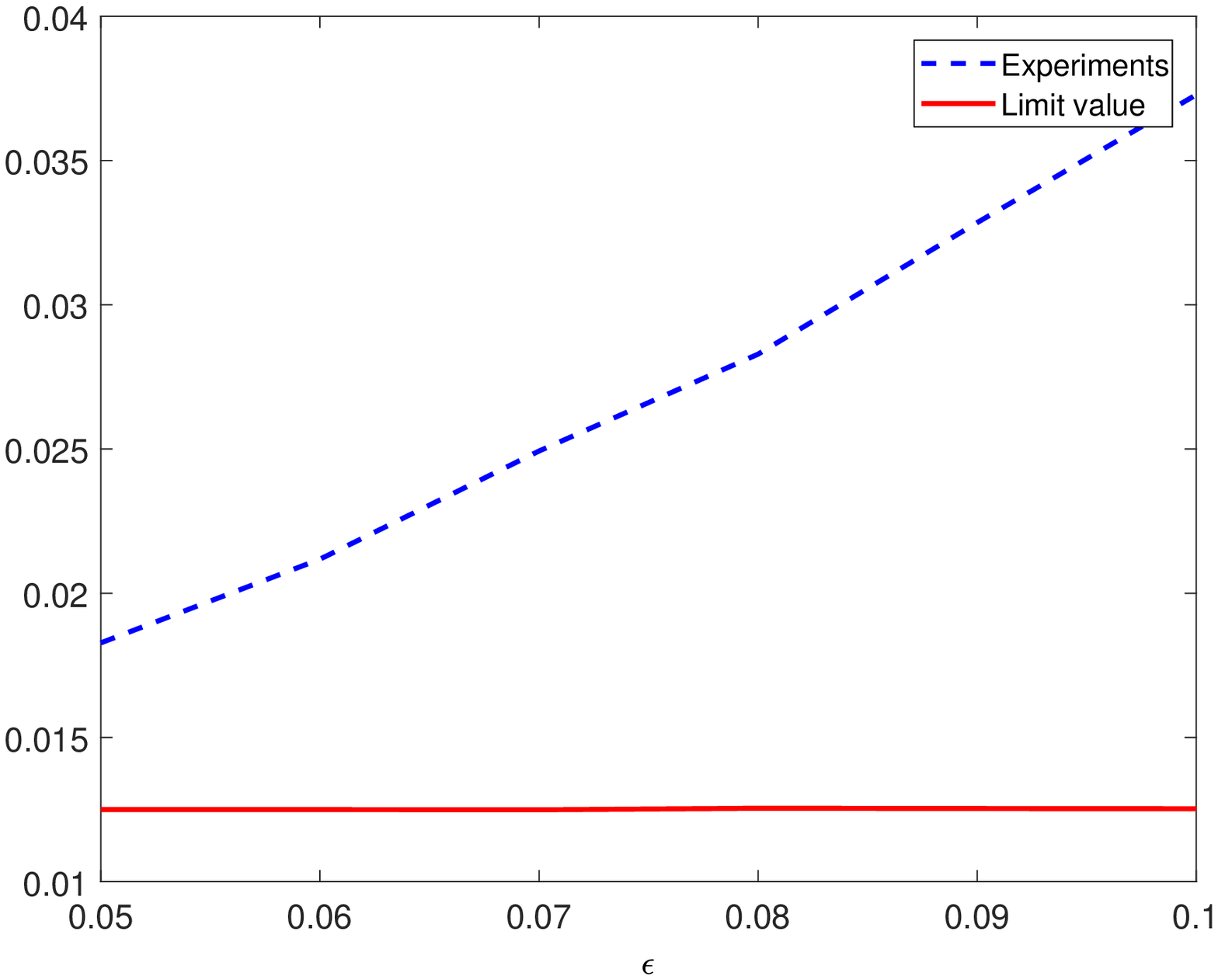}
\caption{Comparison of experiments with asymptotics. The x-axis corresponds to the value of $\epsilon$
and the y-axis to the scaled logarithm of the average exit time.}
\label{fig:Asymptotics}
\end{center}
\end{figure}

At the maxima, the potential achieves the value of 3/8. It is easy to see from the symmetry of the
problem and one-dimensional optimization that the easiest path to reach one maximum from the other
goes through the point $[1/5 \ 3/5 \ 1/5]$ with the value of the potential 7/20.
Thus, in formula (\ref{exitrate}) we have that
$$
\min_{y \in \partial G}(\tilde{F}(\sqrt{x^*}) - \tilde{F}(y)) = 3/8 - 7/20.
$$
In Figure~\ref{fig:Asymptotics} we plot $\epsilon^2\log \bar\tau_{\epsilon}^{exper}$ vs $(3/8 - 7/20)/2$, where $\bar\tau_{\epsilon}^{exper}$ is the average exit time obtained from experiments and $\epsilon$
changes from 0.1 down to 0.05. We note that for values of $\epsilon$ smaller than 0.05, the exit
times become excessively large and we cannot collect enough reliable statistics. Figure~\ref{fig:Asymptotics}
suggests that formula (\ref{exitrate}) describes correctly the asymptotics of the system for
small values of $\epsilon$.

\begin{figure}[h]
\begin{center}
\includegraphics[scale=0.5]{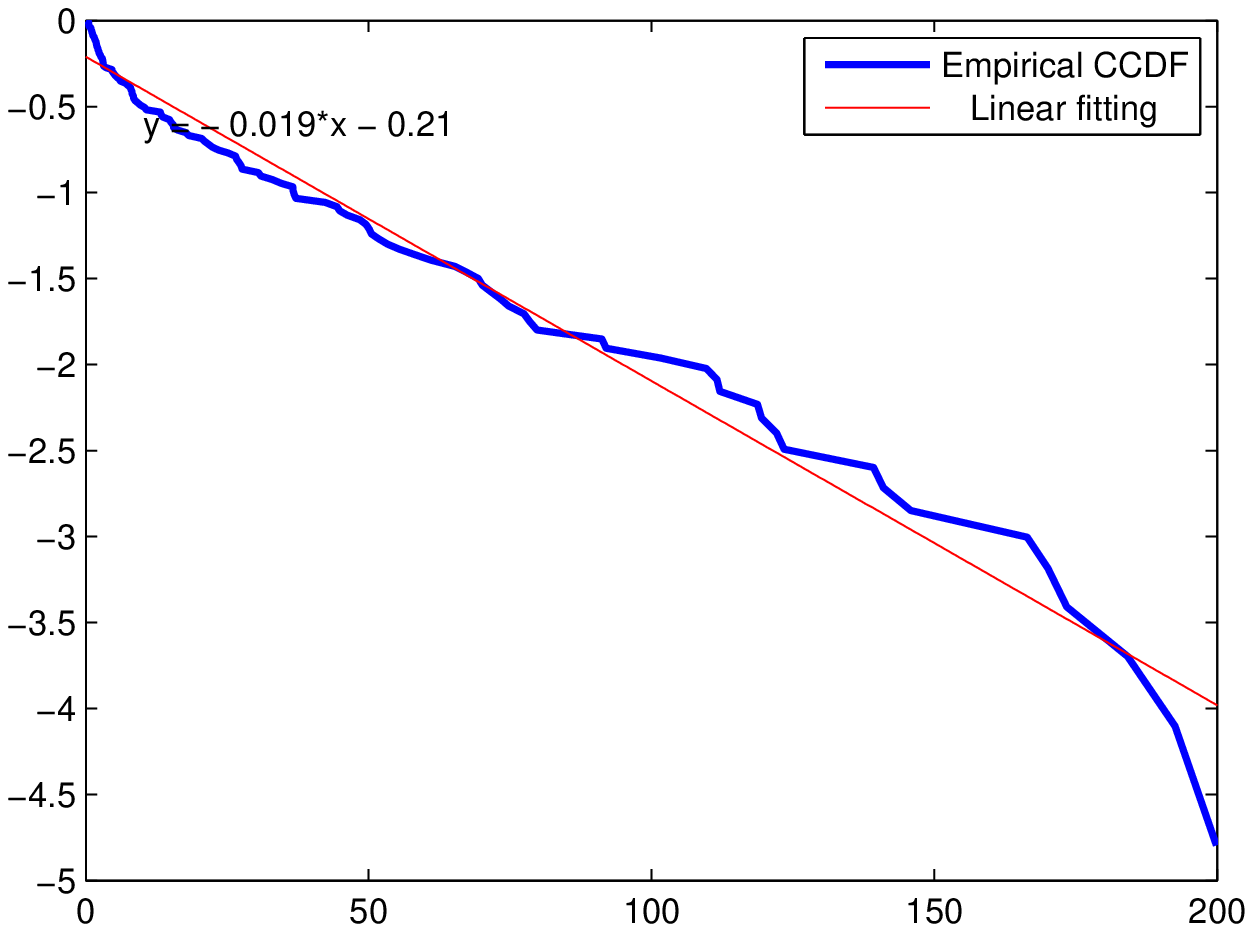}
\caption{Empirical complementary cumulative distribution function of the exit times.
The x-axis corresponds to the value of the exit time.}
\label{fig:A3CCDF}
\end{center}
\end{figure}

We also plot in Figure~\ref{fig:A3CCDF} the empirical complementary cumulative distribution function of the
exit times for $\epsilon=0.1$ in semi-log scale. An approximately linear  shape of the curve is in agreement with
our expectation that the exit times are exponentially distributed. In all numerical experiments of this
section we took the step size $\Delta=0.05$ and the number of steps between one and ten million depending
on the value of $\epsilon$.

\begin{figure}[h]
\begin{center}
\includegraphics[width=\linewidth]{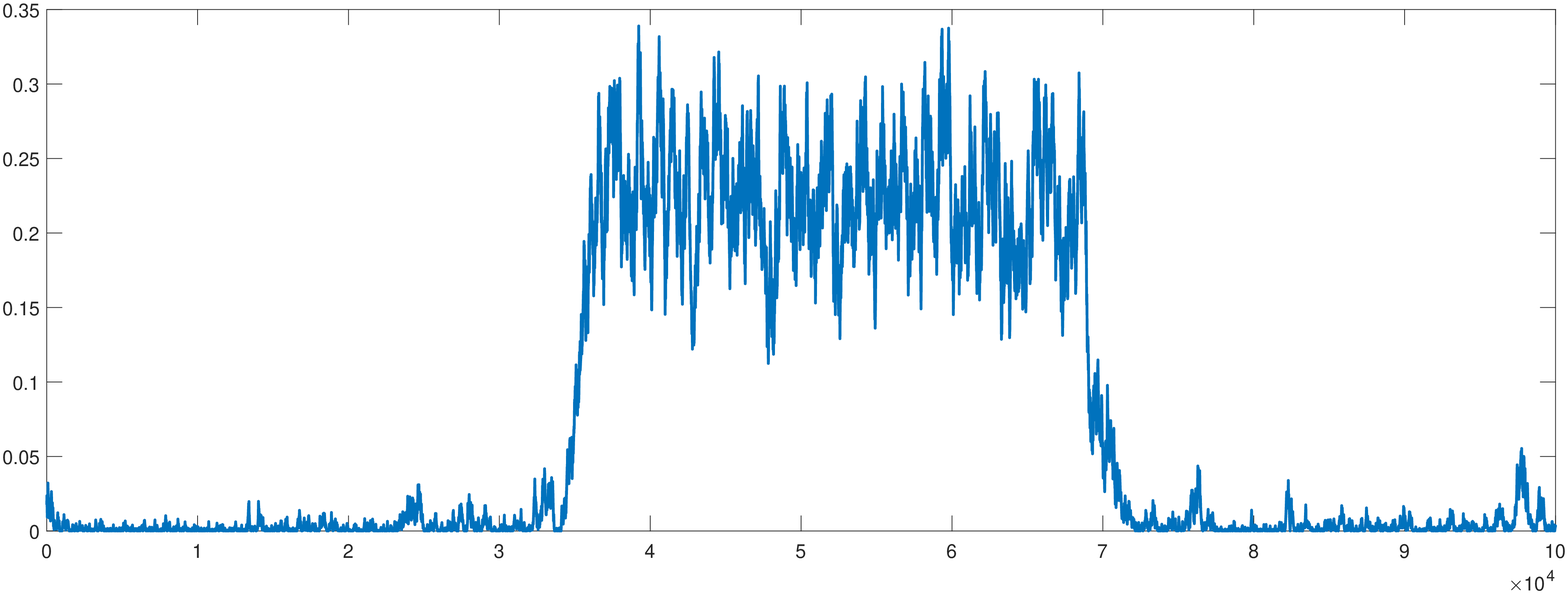}
\caption{A sample of trajectory of one subpopulation in the $G(100,0.25)$ example. The x-axis
corresponds to the time and the y-axis to the subpopulation fraction.}
\label{fig:A100p025TrajBlue}
\end{center}
\end{figure}

Our second numerical example is based on a sample of Erd\H{o}s-R\'{e}nyi random graph with $n=100$ vertices,
with the edge probability $p=0.25$, and with the noise level $\epsilon=0.02$ in the stochastic replicator dynamics.
According to \cite[Chapter~11]{Bollobas}, with high probability, the size of the maximum clique
in Erd\H{o}s-R\'{e}nyi graph is either $\lfloor2\log_{1/p}(n)\rfloor$ or $\lceil2\log_{1/p}(n)\rceil$.
The substitution of $n=100$, $p=0.25$ into the latter formulas indicates that in our example the maximum
clique size is either 6 or 7. Furthermore, using \cite[Eq.~(11.6)]{Bollobas}, we calculate that the
expected number of maximal cliques is 1273 in $G(100,0.25)$. The cliques of sizes 3 and 4 contribute
most to this number. Thus, it is no surprise that most often we observe the cliques of sizes 3 and 4
and rarely 5. In our experiments, we could never observe a clique of size 6 or 7. The system moves from one
metastable state to another and it can take a very long time for the system to come across the metastable state
corresponding to the maximum clique.
A sample of the intermittent evolution of one subpopulation fraction is shown in Figure~\ref{fig:A100p025TrajBlue}.

\begin{figure}[h]
\begin{center}
\includegraphics[width=\linewidth]{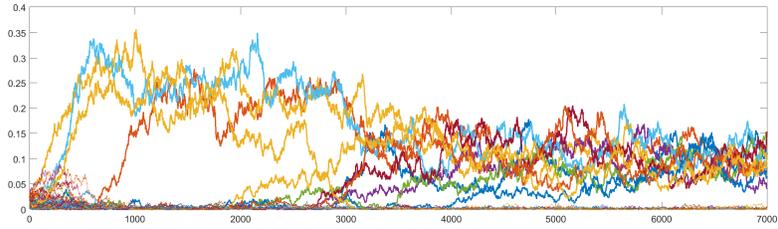}
\caption{Trajectories of subpopulations in the $G(100,0.25)$ example with the planted clique of size 10.
The x-axis corresponds to the time and the y-axis to the subpopulation fractions.}
\label{fig:PlantCliqueTraj}
\end{center}
\end{figure}

\begin{figure}[h]
\begin{center}
\includegraphics[width=\linewidth]{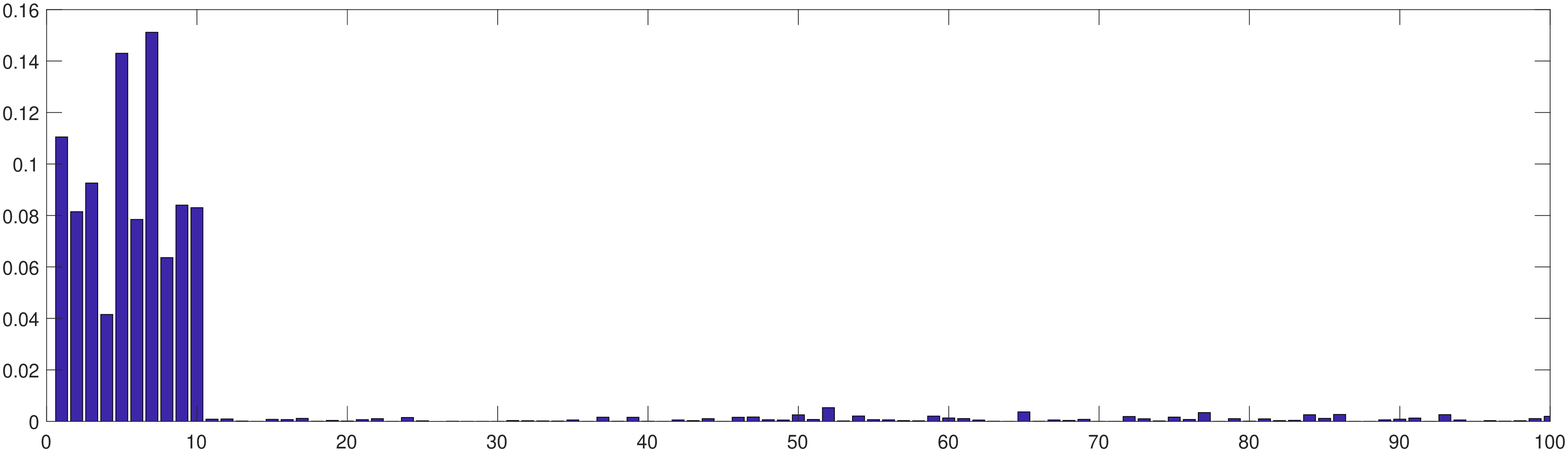}
\caption{The subpopulation fractions at discrete time step 7000 in $G(100,0.25)$ with the planted 
clique of size 10.
The x-axis corresponds to the subpopulation index and the y-axis to the subpopulation fraction. 
The first 10 indices correspond to the planted clique.}
\label{fig:PlantCliqueFreqBar}
\end{center}
\end{figure}

Next, we planted a clique of size 10 in the same realization of $G(100,0.25)$. We have chosen the size 10, because
$\sqrt{n}$ is believed to be the critical scaling for detecting the planted clique \cite{AKS98,Jerrum92}.
It is interesting to observe that in most simulation runs the system state is concentrated on the subpopulations
corresponding to the planted clique
(see as example Figures~\ref{fig:PlantCliqueTraj}~and~\ref{fig:PlantCliqueFreqBar}).

The planted clique experiments, in agreement with the theoretical results from \cite{AKS98,Jerrum92}, indicate that
there is a critical threshold such that, if cardinality of a tightly interacting subpopulation (clique)
exceeds this threshold, the system will fairly quickly converge to the metastable state corresponding to
the planted clique and will remain in that metastable state over long time intervals.


%

\section{Conclusions and future research}

We studied the phenomena of metastability and quasi-invariance as well as intermittent dynamics of a class of stochastic replicator equation in the small noise regime by invoking the Freidlin-Wentzell theory for small noise asymptotics of diffusions and some recent works on quasi-stationarity. We further specialized these results to a noisy replicator dynamics on graphs, with numerical experiments that confirmed our theoretical observations. We see a number of interesting future research directions worth note.
\begin{enumerate}
\item The exit time from the domain of attraction of a stable equilibrium is approximately exponentially distributed. Our example of a stochastic replicator dynamics on graphs shows that an exponential distribution still remains a very good approximation even for the empirical distribution of exit times when the identity of the metastable state is erased, so that each exit episode could be from a different metastable state. In some highly irregular potential functions, however, the exit distribution is observed to be closer to Pareto \cite{Cetal02}. This could be because of clusters of closely spaced metastable states at many length scales, which leads to aggregate effects from a large number of individually exponential random variables. Clearly a more refined analysis is needed to explain this observation.

\item In the context of graph problems, estimating the depth of the deepest well of the potential and finding the precise scaling behavior of the mixing time with the size of the graph (see, e.g., \cite{Vishnoi}) are interesting problems, as is that of characterizing the mixing time in terms of other related stopping times such as hitting times  \cite{Manzo}. Specializing further the graph dynamics to the planted clique problem \cite{AKS98,Jerrum92}, it would be interesting to explore for what scales and after what time the clique becomes detectable by the replicator dynamics.
    These questions have been investigated for a different Metropolis dynamics in \cite{Jerrum92}. It will be interesting
    to carry out a parallel development for our model emphasizing the role of metastability.

\end{enumerate}

\section*{Acknowledgements}

This work was partly supported by EU Project Congas FP7-ICT-2011-8-317672,
CEFIPRA Grants 5100-IT and IFC/DST-Inria-2016-01/448 and ARC Discovery Grant DP 160101236.
We thank the referees and the editor for their careful reading and valuable suggestions which
have immensely improved this article.
We would like to acknowledge the fruitful discussions with Henrik Jensen and Jelena Gruji\'{c}.
We also thank Professors Siva Athreya, K.\ Suresh Kumar, Laurent Miclo, Kavita Ramanan, K.\ S.\ Mallikarjuna Rao
and Alexander Veretennikov for valuable pointers to the literature.

This is the author version of the paper accepted to {\it Dynamic Games and Applications} journal.

%

\end{document}